\documentclass[12pt]{amsart}
\usepackage{amssymb}
\usepackage{amscd}

\newcommand{\Z}{\mathbb{Z}}
\newcommand{\N}{\mathbb{N}}
\newcommand{\C}{\mathbb{C}}

\newcommand{\CalN}{\mathcal{N}}
\newcommand{\CalM}{\mathcal{M}}
\newcommand{\CalV}{\mathcal{V}}

\newcommand{\CalB}{\mathcal{B}}

\newcommand{\CalH}{\mathcal{H}}

\newcommand{\CalO}{\mathcal{O}}
\newcommand{\CalL}{\mathcal{L}}

\newcommand{\gl}{\mathfrak{gl}}

\newcommand{\fraksl}{\mathfrak{sl}}

\newcommand{\End}{\mathrm{End}}

\newcommand{\rk}{\mathrm{rk}\,}

\newcommand{\inv}{\mathrm{inv}}
\newcommand{\Inv}{\mathrm{Inv}}

\newcommand{\leftdiv}{\!\setminus\!}

\newcommand{\Sdt}{\widetilde{S_{d}}}
\newcommand{\Sdh}{\widehat{S_{d}}}
\newcommand{\B}{\CalB}
\newcommand{\Bh}{\widehat{\CalB}}

\newcommand{\m}{\mathbf{m}}

\newcommand{\mmax}{{\mathbf{m}^{\mathrm{max}}}}

\numberwithin{equation}{section}

\newtheorem{theorem}{Theorem}[section]
\newtheorem{corollary}[theorem]{Corollary}

\newtheorem{proposition}[theorem]{Proposition}

\theoremstyle{definition}
\newtheorem*{definition}{Definition}
\newtheorem*{example}{Example}

\newcommand{\bdf}{\begin{definition}}
\newcommand{\edf}{\end{definition}\noindent}
\newcommand{\bex}{\begin{example}}
\newcommand{\eex}{\end{example}\noindent}
\newcommand{\bpr}{\begin{proposition}}
\newcommand{\epr}{\end{proposition}}
\newcommand{\bth}{\begin{theorem}}
\renewcommand{\eth}{\end{theorem}}
\newcommand{\bpf}{\begin{proof}}
\newcommand{\epf}{\end{proof}\noindent}
\newcommand{\bcr}{\begin{corollary}}
\newcommand{\ecr}{\end{corollary}\noindent}
\newcommand{\beq}{\begin{equation}}
\newcommand{\eeq}{\end{equation}}
\newcommand{\bes}{\begin{equation*}}
\newcommand{\ees}{\end{equation*}}
\newcommand{\ben}{\begin{enumerate}}
\newcommand{\een}{\end{enumerate}}
\begin{document}
\title[Nilpotent orbits and Kazhdan-Lusztig polynomials]
{Nilpotent orbits of linear and cyclic quivers 
and Kazhdan-Lusztig poynomials of type A}
\author{Anthony Henderson}
\address{School of Mathematics and Statistics,
University of Sydney, NSW 2006, AUSTRALIA}
\email{anthonyh@maths.usyd.edu.au}
\thanks{This work was supported by Australian Research Council grant DP0344185}
\begin{abstract}
The intersection cohomologies of closures of nilpotent orbits
of linear (respectively, cyclic) quivers are known to be described
by Kazhdan-Lusztig polynomials for the symmetric group (respectively,
the affine symmetric group). We explain how to simplify this
description using a combinatorial cancellation procedure,
and derive some consequences for representation theory.
\end{abstract}
\maketitle
\tableofcontents
\section{Introduction}
This paper is concerned with formulas for the
intersection cohomologies of closures of nilpotent orbits of linear
and cyclic quivers. By fundamental results in geometric representation
theory, these intersection cohomologies control certain features of the
representations of affine Hecke algebras and quantum affine algebras.
There is a well-known formula in the
linear case due to Zelevinsky, using Kazhdan-Lusztig polynomials
of the symmetric group; there is an analogous formula in the
cyclic case due to Lusztig, using Kazhdan-Lusztig polynomials
of the affine symmetric group. The main point of this paper is that
both formulas can be rewritten in terms of Kazhdan-Lusztig polynomials
for different (potentially smaller) symmetric or affine symmetric groups, 
by applying a combinatorial `cancellation' procedure due to Billey
and Warrington. The rewritten formulas in the linear quiver case
have already appeared, in a representation-theoretic guise, in the
work of Suzuki and others; in the cyclic quiver case they are new.

In the remainder of the introduction we will survey the main results
and their representation-theoretic consequences; the other sections
give the proofs, concentrating on the combinatorial side 
of the story. Sections 2 and 4 are purely combinatorial, explaining
the concept of `cancellation' for the symmetric group and affine
symmetric group respectively. Most (perhaps all) of the results in 
Section 2 are known, but we will go over them in detail to provide a
reference for the generalizations to the affine case in Section 4.
Sections 3 and 5 connect these combinatorial results to the problem
of computing intersection cohomology. 

Throughout the paper, all vector spaces, algebras and varieties are over $\C$.

Consider the linear quiver of type $A_{\infty}$, with vertex set $\Z$
and arrows $i\to i+1$ for all $i\in\Z$.
Finite-dimensional representations of this quiver are parametrized by
\emph{multisegments}: a \emph{segment} is a nonempty finite interval
$[i,j]$ in $\Z$, and a multisegment is a finite formal sum of
segments.
Now fix a $\Z$-graded finite-dimensional vector space
$V=\bigoplus_{i\in\Z} V_i$. Let $d_i=\dim V_i$, $d=\dim V$. 
A representation of the quiver on $V$
is simply an element of
\[ \CalN_V=\{\varphi\in\End(V)\,|\,
\varphi(V_{i})\subseteq V_{i+1},\,\forall i\in\Z\}. \]
Two such representations are isomorphic if
they are in the same orbit for $G_V=\{g\in GL(V)\,|\,
g(V_{i})=V_{i},\forall i\}$, acting
on $\CalN_V$ by conjugation. Since all elements of $\CalN_V$ are nilpotent
as endomorphisms of $V$, we call these \emph{nilpotent orbits}.
They are clearly in bijection with $M_{(d_i)}$, the set of multisegments
such that each $i$ occurs $d_i$ times as an element
of a segment. For $\m\in M_{(d_i)}$, let $\CalO_\m$ denote the
corresponding orbit in $\CalN_V$. We put a partial order $\preceq$ on
$M_{(d_i)}$ by setting $\m\preceq\m'$ if and only if $\CalO_\m$ is contained
in the closure $\overline{\CalO_{\m'}}$ of $\CalO_{\m'}$.

The extent to which $\overline{\CalO_{\m'}}$
is singular at the points of $\CalO_{\m}$ is measured by an 
\emph{intersection cohomology} polynomial $IC_{\m,\m'}\in\N[q]$, defined by
\[ IC_{\m,\m'}=\sum_i \dim \CalH_\m^{2i} 
IC(\overline{\CalO_{\m'}})\, q^i, \]
where $IC(\overline{\CalO_{\m'}})$
is the intersection cohomology complex of $\overline{\CalO_{\m'}}$, and
$\CalH_\m^{2i}$ denotes the stalk at a point of $\CalO_\m$
of the $(2i)$th cohomology sheaf (it turns out that all odd-degree
cohomology sheaves of $IC(\overline{\CalO_{\m'}})$ vanish). Note that
$IC_{\m,\m'}$ is nonzero if and only if 
$\m\preceq\m'$, and is $1$ if $\m=\m'$. Hence the
inverse matrix $(IC_{\m,\m'}^{\langle -1\rangle})_{\m,\m'\in M_{(d_i)}}$ of
$(IC_{\m,\m'})_{\m,\m'\in M_{(d_i)}}$ has entries in $\Z[q]$; moreover
$IC_{\m,\m'}^{\langle -1\rangle}$ is zero unless
$\m\preceq\m'$, and is $1$ if $\m=\m'$.

From the viewpoint of geometric representation theory,
the poset $M_{(d_i)}$, together with these IC polynomials, is a model
for certain ``blocks'' of representations of
Lie-theoretic algebras of type $A$. More concretely, the algebras listed
below each have a collection
of finite-dimensional \emph{standard} modules $\{M_\m\,|\,\m\in M_{(d_i)}\}$
and a collection of finite-dimensional simple modules
$\{L_\m\,|\,\m\in M_{(d_i)}\}$, which are related by the following
(equivalent) equations in the Grothendieck group of modules:
\beq \label{multeqn}
\begin{split}
[M_{\m}]&=\sum_{\m'\in M_{(d_i)}}IC_{\m,\m'}(1)\,[L_{\m'}],\ 
\forall\m\in M_{(d_i)},\\
[L_\m]&=\sum_{\m'\in M_{(d_i)}}IC_{\m,\m'}^{\langle -1\rangle}(1)\,[M_{\m'}],
\ \forall\m\in M_{(d_i)}.
\end{split}
\eeq
So the sum of the coefficients of $IC_{\m,\m'}$ is a composition
multiplicity of a standard module; the individual coefficients
record the composition multiplicities in a certain Jantzen-like filtration.
For each algebra, the general definition of standard modules allows
segments of arbitrary complex numbers, not just integers; but the problem
of computing composition multiplicities can be reduced to the integer case.
The algebras in question, and references to the definitions and results,
are as follows.
\ben
\item The affine Hecke algebra $\widehat{\CalH_d}$ attached to
$GL_d$, specialized at a parameter which is not a root of unity
(as in \cite[Definition 12.3.1]{charipressley}). The standard and
simple modules were defined
by Zelevinksy in \cite{zelevinskyhecke}, and \eqref{multeqn} was conjectured
in \cite{zelevinskyconj} (see also \cite{rogawski}).
Ginzburg proved \eqref{multeqn} for standard modules defined
in a geometric way (see \cite[Theorem 8.6.23]{chrissginzburg}). The fact
that Ginzburg's standard modules coincide with Zelevinsky's in the 
Grothendieck group is usually deduced from the Induction Theorem
of Kazhdan and Lusztig (see \cite{ariki} -- nowadays the best version
of the Induction Theorem to use is \cite[Theorem 7.11]{biekii}), 
though it should be regarded as a comparatively easy case of that result.
\item The corresponding degenerate affine Hecke algebra
(as in \cite[Definition 12.3.2]{charipressley}). The definitions
of standard modules (as in \cite[Section 2.2]{suzuki})
are analogous to case (1), and indeed
\eqref{multeqn} in this case can be deduced from case (1) by the results
of Lusztig in \cite{gradedhecke} -- he also gave a proof specific to this
case in \cite{cuspidal}.
\item The quantum affine algebra $U_\epsilon(\widehat{\fraksl_r})$,
specialized at a parameter $\epsilon$ which is 
not a root of unity (as in \cite[Section 12.2A]{charipressley}). 
Here the standard
modules are tensor products of fundamental evaluation modules corresponding
to the segments, so we need $r$ to be greater than or equal to
the length of the longest segment involved (otherwise, we could just stipulate
that any module indexed by a multisegment containing
a segment of length $>r$ is zero). Equation \eqref{multeqn} can be deduced
from case (1) by Frobenius-Schur duality 
(see \cite[Section 12.3D]{charipressley}).
Alternatively, with a geometric definition of standard modules,
\eqref{multeqn} was proved by Ginzburg and Vasserot (see
\cite[Theorem 3]{vasserot}). It then follows that the two kinds
of standard modules are the same in the Grothendieck group (see
\cite[Proposition 18]{vasserot}), which can presumably also be proved
directly.
\item The Yangian $Y(\fraksl_r)$ (as in \cite[Section 12.1A]{charipressley}).
The standard modules as defined by Drinfeld in \cite{drinfeld}
are analogous to those in case (3), and \eqref{multeqn} in this case
can be deduced either from case (3) or case (2) using the results
in \cite{drinfeld}.
\een
This profusion of representation-theoretic meanings of the polynomials
$IC_{\m,\m'}$ and $IC_{\m,\m'}^{\langle -1\rangle}$ is the main 
reason to be interested in computing them; but it is also why, in this
paper, the clear-cut geometric definition is given greater prominence.

A classic result of Zelevinsky (\cite[Corollary 1]{tworemarks},
see Theorem \ref{zelevinskythm} below) 
identifies the polynomials $IC_{\m,\m'}$
with Kazhdan-Lusztig polynomials of the symmetric group $S_d$.
More precisely, it provides an isomorphism of posets
between $M_{(d_i)}$ and a lower ideal 
$M_{(d_i)}'$ of the poset of maximal-length
representatives of $S_{(d_i)}$--$S_{(d_i)}$ double cosets 
(under Bruhat order), where $S_{(d_i)}$ is the parabolic (i.e.\ Young) subgroup
of $S_d$ determined by the composition $d=\sum_i d_i$; and the
polynomials attached to these posets (IC polynomials for $M_{(d_i)}$,
and Kazhdan-Lusztig polynomials for $M_{(d_i)}'$)
coincide under this isomorphism.
Zelevinsky's proof is geometric, embedding
the nilpotent orbits $\CalO_\m$ as open subvarieties of certain
Schubert varieties, and using the fact that the intersection cohomologies
of the latter are described by Kazhdan-Lusztig polynomials; but the
point is that the poset $M_{(d_i)}'$ and its Kazhdan-Lusztig polynomials
can be defined (and, in principle, computed) purely combinatorially.

More recent work of Suzuki (\cite{suzuki}) implicitly generalizes
this result, providing a family of
poset isomorphisms between various upper ideals of $M_{(d_i)}$ and 
combinatorially-defined posets. 
To explain this we adopt the notation of \cite{orellanaram},
which views multisegments as ``generalized skew-shapes''.
For $\lambda,\mu\in\Z^k$, write $\lambda\supseteq\mu$ if $\lambda_i\geq\mu_i$
for all $1\leq i\leq k$, and if this holds define a multisegment
\beq \label{skewshapeeqn}
\lambda/\mu=\sum_{i=1}^k\, [\mu_i-i+1,\lambda_i-i], 
\eeq
where any ``empty segments'' of the form $[s+1,s]$ are ignored.
The reason for the notation is that if $\lambda$ and $\mu$ are 
\emph{partitions},
i.e.\ $\lambda_1\geq\cdots\geq\lambda_k\geq 0$ and similarly for $\mu$,
then the segments are exactly the rows of the skew-shape diagram usually
called $\lambda/\mu$, where each box is replaced by its \emph{content}.
Since the order of
terms in \eqref{skewshapeeqn} is unimportant,
$\lambda/\mu=(w\cdot\lambda)/
(w\cdot\mu)$ for all $w\in S_k$, where the ``dot action'' of $S_k$ on $\Z^k$
is defined as usual by
\beq 
(w\cdot\lambda)_i-i=\lambda_{w^{-1}(i)}-w^{-1}(i).
\eeq
A fundamental domain for this dot action is
\[ D_k:=\{\lambda\in\Z^k\,|\,
\lambda_1-1\geq\lambda_2-2\geq\cdots\geq\lambda_k-k\}, \]
so we can write any multisegment
in the \emph{standard form}
\beq \label{standardformeqn}
\lambda/(w\cdot\mu) \text{ where }\left\{\begin{array}{c}
\lambda,\mu\in D_k,\ w\in S_k,\ \lambda\supseteq w\cdot\mu,\text{ and}\\
w\text{ has maximal length in }W_\lambda w W_\mu.\end{array}\right.  
\eeq
Here $W_\lambda$ and $W_\mu$ are the stabilizers of $\lambda$ and $\mu$
for the dot action, which are clearly parabolic subgroups of $S_k$.
Note that this expression in standard form is not uniquely
determined by the multisegment, but rather by the multisegment together
with a chosen multiset of empty segments.
\bex
Let $\m$ be the multisegment $[1,2]+[2,2]+[3,3]$. The most economical way
to express this in the form \eqref{standardformeqn} is to take $k=3$,
$\lambda=(4,4,5)$, $\mu=(3,3,3)$, and $w$ to be the transposition $(2,3)$.
Another way is to take $k=4$, $\lambda=(4,4,5,5)$, $\mu=(3,3,4,4)$,
and $w$ to be the transposition $(2,4)$; this effectively adds the empty 
segment $[2,1]$.
\eex
If $\lambda,\mu\in D_k$, we define 
\bes
\begin{split} 
S_k[\lambda,\mu]&=\{w\in S_k\,|\,\lambda\supseteq w\cdot\mu\}\text{ and}\\
S_k[\lambda,\mu]^\circ&=\{w\in S_k[\lambda,\mu]\,|\,w
\text{ has maximal length in }W_\lambda w W_\mu\}.
\end{split}
\ees
These are posets under Bruhat order; in fact we will see that
$S_k[\lambda,\mu]$ is a lower ideal of $S_k$, so in particular
$S_k[\lambda,\mu]\neq\emptyset\Leftrightarrow\lambda\supseteq\mu$.
With this notation, the generalized form of Zelevinsky's result can
be stated as follows.
\bth \label{introfinmainthm}
Let $\lambda,\mu\in D_k$ be such that 
$\lambda\supseteq\mu$, $\lambda/\mu\in M_{(d_i)}$.
\ben
\item The map $w\mapsto \lambda/(w\cdot\mu)$ is an isomorphism of posets
between $S_k[\lambda,\mu]^\circ$ and 
$\{\m'\in M_{(d_i)}\,|\,\lambda/\mu\preceq\m'\}$.
\item For $w,w'\in S_k[\lambda,\mu]^\circ$,
$IC_{\lambda/(w\cdot\mu),\lambda/(w'\cdot\mu)}=P_{w,w'}$,
a Kazhdan-Lusztig polynomial of $S_k$.
\een
\eth
\noindent
Zelevinsky's original result is the special case where 
$k=d$, $\lambda$ is such that each integer $i$
occurs $d_i$ times in $(\lambda_1-1,\cdots,\lambda_d-d)$, and
$\mu=\lambda-(1,1,\cdots,1)$. In this case $\lambda/\mu$ is the trivial
multisegment $\sum_i d_i [i,i]$ (corresponding to the zero orbit),
so the image of the isomorphism in part (1) is all of $M_{(d_i)}$;
the parabolic subgroups $W_\lambda$ and $W_\mu$
both equal $S_{(d_i)}$, and $S_d[\lambda,\mu]^\circ$ is the poset $M_{(d_i)}'$ 
mentioned above.

The fact that Theorem \ref{introfinmainthm} is true in general means
that in expressing $\m$ and $\m'$ as $\lambda/(w\cdot\mu)$ and
$\lambda/(w'\cdot\mu)$, any of  
the empty segments which occur in the Zelevinsky case can
be ``cancelled'' without changing the Kazhdan-Lusztig polynomial
(or, indeed, new ones can be added). 
For example, at the extreme, Theorem \ref{introfinmainthm}
shows that $IC_{\m,\m'}$ can be identified with a Kazhdan-Lusztig
polynomial of $S_{k(\m)}$, where $k(\m)$ is the number of segments
of $\m$ (smaller than $d$, unless $\m$ is trivial).
In Section 3 we will
use a result of Billey and Warrington, which provides for just
such cancellations in Kazhdan-Lusztig polynomials of symmetric groups,
to deduce Theorem \ref{introfinmainthm} from
Zelevinsky's theorem. The essence of the result is stated in Theorem 
\ref{finmainthm}, using some matrix notation which will be introduced in \S2.

As already mentioned, Theorem \ref{introfinmainthm} cannot be considered new,
since in the context of the representation theory of the degenerate affine
Hecke algebra (case (2) above) it follows from Suzuki's results in 
\cite{suzuki}. With notation as in the Theorem, he defines an 
exact functor $F_\lambda$ from
the category $\CalO$ of representations of $\gl_k$ to the category
of finite-dimensional modules for the degenerate affine Hecke algebra
associated to $GL_d$, and shows that it takes the Verma module
$M(w\cdot\mu)$ to the standard module $M_{\lambda/(w\cdot\mu)}$
and the simple module $L(w\cdot\mu)$ to the simple module
$L_{\lambda/(w\cdot\mu)}$ for all $w\in S_k[\lambda,\mu]^\circ$.
Part (1) and the $q=1$ specialization of part (2) of Theorem 
\ref{introfinmainthm} then follow from the known
Kazhdan-Lusztig conjecture for $\gl_k$ (combine 
\cite[(5.2.1) and (5.2.2)]{suzuki} -- the historical remarks following
\cite[(5.2.3)]{suzuki} properly apply only to the Zelevinsky case).
Moreover, by \cite[Theorem 5.3.5]{suzuki} the Kazhdan-Lusztig
polynomials $P_{w,w'}$ for $w'\in S_k[\lambda,\mu]^\circ$ record
multiplicities in a Jantzen-type filtration of $M_{\lambda/(w\cdot\mu)}$,
whose definition clearly depends only on the multisegment
(i.e.\ not on the empty segments); since (2) of
Theorem \ref{introfinmainthm} is true in the Zelevinsky case, it must
be true in general. (As well as \cite{suzuki}, see \cite{orellanaram}
and \cite{arakawa}
for the analogous results in the case of the affine Hecke algebra
and Yangian respectively).

One corollary concerns those
multisegments $\lambda/\mu$ where
$\lambda,\mu\in D_k$ satisfy $W_\lambda=W_\mu=\{1\}$, i.e.
$\lambda_1\geq\lambda_2\geq\cdots\geq\lambda_k$,
$\mu_1\geq\mu_2\geq\cdots\geq\mu_k$ (these are the ``placed skew-shapes''
of \cite{ram}): for such $\lambda/\mu$, it follows from
Theorem \ref{introfinmainthm} that
\beq \label{signinveqn}
IC_{\lambda/\mu,\lambda/(w\cdot\mu)}^{\langle -1\rangle}=\varepsilon(w),
\text{ for all }w\in S_k[\lambda,\mu],
\eeq
where $\varepsilon$ denotes the sign character.
Thus the corresponding simple modules (called \emph{calibrated}
for the affine Hecke algebra in \cite{ram} and \emph{tame}
for the Yangian in \cite{nazarovtarasov}) can be written as an
alternating sum of standard modules in the Grothendieck group.
Representation-theoretically, this reflects the existence of a
BGG-like resolution of these simple modules (transferred by
the appropriate functor from the BGG resolution of the $\gl_k$-module
$L(\mu)$) -- see \cite[Theorem 5.1.1]{suzuki} and
\cite[(4.13)]{orellanaram}.

The justification for re-proving Theorem \ref{introfinmainthm}
in \S3 below is that the combinatorics involved
generalizes immediately to the case of cyclic quivers,
as we will now explain.

Fix a positive integer $n$, and
consider the cyclic quiver of type $\widetilde{A_{n-1}}$, with vertex set
$\Z/n\Z$ and arrows $\bar i\to\overline{i+1}$ for all $\bar i\in\Z/n\Z$.
Finite-dimensional \emph{nilpotent} representations
of this quiver are parametrized by multisegments as before, except that
there is no difference between segments $[i,j]$ and $[i',j']$ when
$i'-i=j'-j$ is a multiple of $n$.
Fix a $(\Z/n\Z)$-graded finite-dimensional
vector space $V=\bigoplus_{\bar i\in\Z/n\Z} V_{\bar i}$, and set
$d_{i}=\dim V_{\bar i}$, $d=\dim V$. We define
\[ \CalN_V=\{\varphi\in\End(V)\,|\,
\varphi(V_{\bar{i}})\subseteq V_{\overline{i+1}},\,\forall\bar{i}\in\Z/n\Z,\
\varphi\text{ nilpotent}\}, \]
and consider $G_V$-orbits in $\CalN_V$. These are in bijection with
$M_{(d_i),n}$, the set of multisegments 
(in this modulo $n$ sense) such that each congruence class $\bar i$
occurs $d_i$ times among the elements of the segments.
For $\m\in M_{(d_i),n}$, let $\CalO_\m$ denote the corresponding nilpotent
orbit, and define a partial order $\preceq$ and
polynomials $IC_{\m,\m'}\in\N[q]$ and
$IC_{\m,\m'}^{\langle -1\rangle}\in\Z[q]$ in the same way
as before.

These polynomials too have representation-theoretic significance.
The specialized quantum affine algebra $U_\zeta(\widehat{\fraksl_r})$, where
$\zeta^2$ is a primitive $n$th root of $1$, has a collection
of standard modules $\{M_{\m}\,|\,\m\in M_{(d_i),n}\}$ and a collection
of simple modules
$\{L_{\m}\,|\,\m\in M_{(d_i),n}\}$,
satisfying the equivalent equations
\beq \label{affmulteqn}
\begin{split}
[M_{\m}]&=\sum_{\m'\in M_{(d_i),n}}IC_{\m,\m'}(1)\,[L_{\m'}],\ 
\forall\m\in M_{(d_i),n},\\
[L_\m]&=\sum_{\m'\in M_{(d_i),n}}IC_{\m,\m'}^{\langle -1\rangle}(1)\,[M_{\m'}],
\ \forall\m\in M_{(d_i),n}.
\end{split}
\eeq
(See \cite[Theorem 3]{vasserot} -- again, for small $r$ 
we have to disregard multisegments containing
a segment of length $>r$.) The same is true for the affine Hecke algebra
$\widehat{\CalH_d}$ specialized at a primitive $n$th root of unity,
except that there the simple modules are parametrized by the smaller set
of \emph{aperiodic} multisegments
(see \cite[Section 2]{ltv}), so we have to set $[L_\m]=0$ if $\m$ is not
aperiodic.

The analogue of Zelevinsky's result for cyclic quivers was proved by
Lusztig in \cite[\S11]{quivers1} (it is stated below as Theorem 
\ref{lusztigthm}). This identifies $IC_{\m,\m'}$ with a Kazhdan-Lusztig
polynomial of the affine symmetric group $\widetilde{S_d}$
(the Coxeter group of type $\widetilde{A_{d-1}}$). In Section 4,
we will show that a version of Billey and Warrington's cancellation works for
the affine symmetric group. As a consequence, we get an analogue
of Theorem \ref{introfinmainthm} in this setting, for which a
representation-theoretic proof does not yet exist.

To state it requires extending the dot action of $S_k$ on $\Z^k$
to $\widetilde{S_k}$, so that the extra Coxeter generator $s_0$ acts by
\[ (s_0\cdot\lambda)_i=\left\{\begin{array}{cl}
\lambda_k-k+1+n,&\text{ if $i=1$,}\\
\lambda_i,&\text{ if $2\leq i\leq k-1$,}\\
\lambda_1-1+k-n,&\text{ if $i=k$.}
\end{array}\right. \]
It is then clear that $(w\cdot\lambda)/(w\cdot\mu)=\lambda/\mu$
for all $w\in\widetilde{S_k}$, where
the multisegments are now interpreted in the modulo $n$ sense.
A fundamental domain for the action of $\widetilde{S_k}$ on $\Z^k$ is
\[ \widetilde{D_k}:=\{\lambda\in\Z^k\,|\,
\lambda_1-1\geq\lambda_2-2\geq\cdots\geq\lambda_k-k\geq\lambda_1-n-1\}, \]
and the corresponding standard form for multisegments is
\beq \label{affstandardformeqn}
\lambda/(w\cdot\mu) \text{ where }\left\{\begin{array}{c}
\lambda,\mu\in \widetilde{D_k},\ w\in \widetilde{S_k},\ 
\lambda\supseteq w\cdot\mu,\text{ and}\\
w\text{ has maximal length in }\widetilde{W_\lambda} w \widetilde{W_\mu},
\end{array}\right.  
\eeq
where $\widetilde{W_\lambda}$ and $\widetilde{W_\mu}$ denote the stabilizers
of $\lambda$ and $\mu$ in $\widetilde{S_k}$ (proper parabolic subgroups,
hence finite).
For $\lambda,\mu\in\widetilde{D_k}$, we define
\bes
\begin{split} 
\widetilde{S_k}[\lambda,\mu]&=
\{w\in \widetilde{S_k}\,|\,\lambda\supseteq w\cdot\mu\}\text{ and}\\
\widetilde{S_k}[\lambda,\mu]^\circ&=\{w\in \widetilde{S_k}[\lambda,\mu]\,|\,w
\text{ has maximal length in }\widetilde{W_\lambda} w \widetilde{W_\mu}\}.
\end{split}
\ees
We will see in \S5 that, as in the symmetric group case, 
$\widetilde{S_k}[\lambda,\mu]$ is a (finite) 
lower ideal of $\widetilde{S_k}$ for
Bruhat order. We can now state a generalization of Lusztig's result.
\bth \label{introaffmainthm}
Let $\lambda,\mu\in \widetilde{D_k}$ be such that $\lambda\supseteq\mu$,
$\lambda/\mu\in M_{(d_i),n}$.
\ben
\item The map $w\mapsto \lambda/(w\cdot\mu)$ is an isomorphism of posets
between $\widetilde{S_k}[\lambda,\mu]^\circ$ and 
$\{\m'\in M_{(d_i),n}\,|\,\lambda/\mu\preceq\m'\}$.
\item For $w,w'\in \widetilde{S_k}[\lambda,\mu]^\circ$,
$IC_{\lambda/(w\cdot\mu),\lambda/(w'\cdot\mu)}=P_{w,w'}$,
a Kazhdan-Lusztig polynomial of $\widetilde{S_k}$.
\een
\eth
\noindent
(An alternative statement using matrix notation is given in Theorem
\ref{mainthm}.)
It is natural to wonder whether there is a representation-theoretic functor
which ``explains'' this Theorem too.

As in the linear quiver case, Theorem \ref{introaffmainthm} implies
that $IC_{\m,\m'}$ can be identified with a Kazhdan-Lusztig polynomial
of $\widetilde{S_{k(\m)}}$, where $k(\m)$ is the number of segments
of $\m$; this immediately implies the main result of \cite{mytworow},
that $IC_{\m,\m'}=1$ when $\m\preceq\m'$, $k(\m)=2$.

Another consequence of Theorem \ref{introaffmainthm}
is an analogue of \eqref{signinveqn}, concerning those multisegments 
$\lambda/\mu$ where $\lambda,\mu\in\widetilde{D_k}$
satisfy $\widetilde{W_\lambda}=\widetilde{W_\mu}=\{1\}$; this means that
$\lambda_1\geq\lambda_2\geq\cdots\lambda_k\geq\lambda_1-n+k$,
$\mu_1\geq\mu_2\geq\cdots\geq\mu_k\geq\mu_1-n+k$. For such $\lambda/\mu$,
it follows from Theorem \ref{introaffmainthm} that
\beq \label{affsigninveqn}
IC_{\lambda/\mu,\lambda/(w.\mu)}^{\langle -1\rangle}=\varepsilon(w),
\text{ for all }w\in \widetilde{S_k}[\lambda,\mu].
\eeq
So once more the corresponding simple modules can be written as an
alternating sum of standard modules in the Grothendieck group;
probably this indicates a BGG-like resolution.

Theorems \ref{introfinmainthm} and \ref{introaffmainthm} combine well
with the method used by Varagnolo and Vasserot in \cite{vv} to determine
the decomposition numbers
of $U_\zeta(\gl_r)$ where $\zeta^2$ is a primitive $n$th root of $1$.
Suppose we want to compute the multiplicity of the simple module
$L_\zeta(\mu')$ in the Weyl module $V_\zeta(\lambda')$, where
$\lambda$ and $\mu$ are partitions with at most $k$ parts all
of size $\leq r$, and
$\lambda'$ and $\mu'$ are the transpose partitions (regarded as
dominant integral weights for $\gl_r$). By definition, $V_\zeta(\lambda')$
is the specialization at $\zeta$ of the simple module $V_q(\lambda')$
for the generic $U_q(\gl_r)$. Now using a suitable normalization of the
evaluation map $U_q(\widehat{\fraksl_r})\to U_q(\gl_r)$, we can regard
$V_q(\lambda')$ as the simple $U_q(\widehat{\fraksl_r})$-module
$L_{\lambda/0}$ (see \cite[Section 12.2]{vv}). 
By \eqref{signinveqn}, we have the equation
\beq 
[L_{\lambda/0}]=\sum_{w\in S_k[\lambda,0]}\varepsilon(w)\,
[M_{\lambda/(w\cdot 0)}].
\eeq
Now let $w_\lambda, w_\mu, w_0\in\widetilde{S_k}$ be such that
$w_\lambda\cdot\lambda,w_\mu\cdot\mu,w_0\cdot 0\in\widetilde{D_k}$.
As noted in \cite[Section 12.3]{vv}, the specialization at 
$\zeta$ of the standard module
$M_{\lambda/(w\cdot 0)}$ is merely the $U_\zeta(\widehat{\fraksl_r})$-standard 
module of the same name, which in standard form is
$M_{(w_\lambda\cdot\lambda)/(w_\lambda w w_0^{-1})^{\circ}\cdot(w_0\cdot 0)}$,
where $(w_\lambda w w_0^{-1})^{\circ}$ is the longest element of
$\widetilde{W_{w_\lambda\cdot\lambda}}w_\lambda w w_0^{-1}
\widetilde{W_{w_0\cdot 0}}$.
So in the Grothendieck group of $U_\zeta(\widehat{\fraksl_r})$-modules,
\beq
[V_\zeta(\lambda')]=\sum_{w\in S_k[\lambda,0]}\varepsilon(w)\,
[M_{(w_\lambda\cdot\lambda)/(w_\lambda w w_0^{-1})^{\circ}\cdot(w_0\cdot 0)}].
\eeq
Now as noted in \cite[Section 12.2]{vv}, $L_\zeta(\mu')$ regarded as
a simple $U_\zeta(\widehat{\fraksl_r})$-module is $L_{\mu/0}=
L_{(w_\mu\cdot\mu)/(w_\mu w_0^{-1})^{\circ}\cdot(w_0\cdot 0)}$, where
$(w_\mu w_0^{-1})^{\circ}$ is the longest element of
$\widetilde{W_{w_\mu\cdot\mu}}w_\mu w_0^{-1}\widetilde{W_{w_0\cdot 0}}$. 
Using Theorem \ref{introaffmainthm}, we obtain
\beq \label{bigeqn}
\begin{split}
[V_\zeta(\lambda'):L_\zeta(\mu')]
=\left\{\begin{array}{cl}
\sum_{w\in S_k}\varepsilon(w)\,P_{(w_\lambda w w_0^{-1})^{\circ},
(w_\mu w_0^{-1})^{\circ}}(1),&\text{ if $\mu\in\widetilde{S_k}\cdot\lambda$,}\\
0,&\text{ otherwise.}
\end{array}\right.
\end{split}
\eeq
In the first case, summing over all of $S_k$ 
rather than just $S_k[\lambda,0]$ introduces
no new terms, since $(w_\mu w_0^{-1})^{\circ}\in\widetilde{S_k}[w_\mu\cdot\mu,
w_0\cdot 0]=\widetilde{S_k}[w_\lambda\cdot\lambda,w_0\cdot 0]$, 
so the Kazhdan-Lusztig polynomial can only be nonzero when
$w_\lambda w w_0^{-1}\in\widetilde{S_k}[w_\lambda\cdot\lambda,w_0\cdot 0]$, 
i.e.\ $\lambda\supseteq w\cdot 0$.

In the special case that $\lambda$ and $\mu$ have trivial stabilizers in 
$\widetilde{S_k}$ (i.e.\ $\lambda_1-1,\cdots,\lambda_k-k$ have different 
residues modulo $n$, and similarly for $\mu$ -- 
this requires $k\leq n$, which automatically implies
$w_0=1$), \eqref{bigeqn} becomes
\beq
\begin{split}
[V_\zeta(\lambda'):L_\zeta(\mu')]
=\left\{\begin{array}{cl}
\sum_{w\in S_k}\varepsilon(w)\,P_{w_\lambda w,
w_\mu}(1),&\text{ if $\mu\in\widetilde{S_k}\cdot\lambda$,}\\
0,&\text{ otherwise.}
\end{array}\right.
\end{split}
\eeq
This is the form of the answer given by Soergel
in \cite[Conjecture 7.1]{soergel}
for the equivalent
problem of computing tilting module multiplicities for $U_\zeta(\gl_k)$.

\noindent
\textit{Acknowledgements. }
This work grew out of stimulating conversations with my colleagues
A.~Mathas and A.~Molev, and I would like to thank them for their help.
I am also grateful to G.~Lusztig, K.~McGerty, A.~Parker,
A.~Ram, J.~Ramagge, and G.~Warrington for their valuable comments on an earlier
version of this paper.
\section{Cancellation for the symmetric group}
In this section we explain the combinatorial result of Billey and
Warrington on which our approach depends.
Fix a positive integer $d$, and let $S_d$ be the group of permutations
of $[1,d]=\{1,\cdots,d\}$. For $i\in[1,d-1]$, we define $s_i\in S_d$
to be the transposition interchanging $i$ and $i+1$; as everyone knows,
$s_1,\cdots,s_{d-1}$ form a set of Coxeter generators for $S_d$ of
type $A_{d-1}$. We thus have a length function $\ell:S_d\to\N$,
a Bruhat order $\leq$, and Kazhdan-Lusztig polynomials $P_{y,w}\in\N[q]$
for $y,w\in S_d$ (which are nonzero iff $y\leq w$). Good references
for Kazhdan-Lusztig polynomials are \cite[Chapter 7]{humphreys}
and \cite{soergel} (where the notation is somewhat different).

The length function and the Bruhat order have well-known 
combinatorial descriptions.
Define the \emph{inversion statistics} 
\[ \inv_i(w)=|\{i'<i\,|\,w(i')>w(i)\}|,\
\Inv_i(w)=|\{i'>i\,|\,w(i')<w(i)\}|, \]
for any $w\in S_d$ and $i\in[1,d]$.
These are related by $\Inv_i(w)=\inv_i(w)+w(i)-i$. Then
\beq \label{finlengtheqn}
\ell(w)=\sum_{i\in [1,d]}\inv_i(w)=\sum_{i\in [1,d]}\Inv_i(w).
\eeq
A special case of Bruhat order is that for all $i\in [1,d-1]$,
\beq 
ws_i<w\text{ if and only if }
w(i)>w(i+1).
\eeq 
The general description, due to Deodhar, is as follows:
\bpr \label{finbruhatprop}
If $y,w\in S_d$, $y\leq w$ if and only if for all $i,j\in[1,d]$,
\[ |\{i'\leq i\,|\,y(i')\geq j\}|\leq |\{i'\leq i\,|\,w(i')\geq j\}|. \]
\epr
\noindent
In other words, for all $i\in [1,d]$ and $m\in [1,i]$, 
the $m$th largest element in $y[1,i]$
is less than or equal to the $m$th largest element in $w[1,i]$.
If $y\leq w$, we write $[y,w]$ for the Bruhat interval
$\{x\in S_d\,|\,y\leq x\leq w\}$.

We now come to the key definition.
\bdf
If $y\leq w$ in $S_d$, we say that $i\in[1,d]$ is 
\emph{cancellable} for
the interval $[y,w]$ if $y(i)=w(i)$, $\inv_{i}(y)=\inv_{i}(w)$, and
$\Inv_i(y)=\Inv_i(w)$. (Clearly any two of these conditions imply the third.)
\edf
The reason for the name `cancellable' is that Bruhat order and
Kazhdan-Lusztig polynomials are preserved under the operation of
`cancelling the common action on $i$' from 
the permutations
in question, in the following sense. For all $i\in[1,d]$, let
$\sigma_{i}:[1,d]\setminus\{i\}\to[1,d-1]$ be
the unique order-preserving bijection.
For $w\in S_d$, we define $w^{\hat{i}}\in S_{d-1}$ by
\[ w^{\hat{i}}=\sigma_{w(i)}\circ w\circ\sigma_{i}^{-1}. \]
It is clear from either formula in \eqref{finlengtheqn} that
\beq \label{lengthdropeqn}
\ell(w^{\hat{i}})=\ell(w)-\inv_i(w)-\Inv_i(w). 
\eeq
The following result combines Lemmas 17 and 39 of \cite{billeywarrington}, but
we will spell out the proof for later reference.
\bpr \label{fincancprop}
Suppose that $i$ is cancellable for $[y,w]$.
\ben
\item For any $x\in [y,w]$, $x(i)=y(i)$ and $\inv_{i}(x)=\inv_{i}(y)$.
Hence $i$ is cancellable for any sub-interval of $[y,w]$.
\item $x\mapsto x^{\hat{i}}$ is an isomorphism of posets between $[y,w]$
and $[y^{\hat{i}},w^{\hat{i}}]$, which reduces all lengths by the same
amount.
\item For any $u,v\in [y,w]$, $P_{u,v}=P_{u^{\hat{i}},v^{\hat{i}}}$.
\een
\epr
\bpf
Set $j=y(i)=w(i)$, $m=\inv_{i}(y)+1=\inv_{i}(w)+1$, and 
suppose $y\leq x\leq w$. 
Now $y[1,i]$ and $w[1,i]$ each have exactly
$m$ elements $\geq j$ and $m-1$ elements $>j$. 
By Proposition \ref{finbruhatprop}, the same is true of
$x[1,i]$. Similarly,
$y[1,i-1]$ and $w[1,i-1]$ each have exactly
$m-1$ elements $\geq j$ and $m-1$ elements $>j$,
so the same is true of $x[1,i-1]$. Thus $x(i)=j$ and $\inv_i(x)=m-1$,
proving (1).
Moreover, it is clear from Proposition \ref{finbruhatprop} that
$y^{\hat{i}}\leq x^{\hat{i}}\leq w^{\hat{i}}$.
The construction of the map $[y^{\hat{i}},w^{\hat{i}}]\to
[y,w]:x\mapsto\tilde{x}$ inverse to $x\mapsto x^{\hat{i}}$ is easy:
\[ \tilde{x}(i')=\left\{\begin{array}{cl}
j,&\text{ if $i'=i$,}\\
\sigma_{j}^{-1}(x(\sigma_{i}(i'))),
&\text{ if $i'\neq i$.}
\end{array}\right. \]
This proves the isomorphism part of (2), and the statement about lengths
follows from (1).
In light of parts (1) and (2), it clearly suffices to prove
(3) in the case $u=y$, $v=w$.
We prove this by induction on $\ell(w)$, it being trivial if $w=1$.
Choose one of the Coxeter generators, say $s$, such that $ws<w$. 
We now have three cases.\newline
\textbf{Case 1: }$s=s_{i-1}$. This means that $w(i-1)>j$, so
$w[1,i-2]$ has only $m-2$ elements $> j$.
Therefore the same is true of $y[1,i-2]$, so $y(i-1)>j$,
i.e.\ $ys<y$. Moreover, $ws[1,i-1]$ has only
$m-2$ elements $>j$, so $y\not\leq ws$.
Under these circumstances we have (see \eqref{kleqn} below)
\beq P_{y,w}=P_{ys,ws}.
\eeq
Obviously $i-1$ is cancellable for $[ys,ws]$, so by the
induction hypothesis, $P_{ys,ws}=
P_{(ys)^{\widehat{i-1}},(ws)^{\widehat{i-1}}}$.
But $(ys)^{\widehat{i-1}}=y^{\hat{i}}$ and
$(ws_{i-1})^{\widehat{i-1}}=w^{\hat{i}}$, so we have the result.\newline
\textbf{Case 2: }$s=s_i$. This means that $w(i+1)<j$, so
$w[1,i+1]$ has only $m-1$ elements $>j$.
Therefore the same is true of $y[1,i+1]$, so $y(i+1)<j$, i.e.\ $ys<y$.
Moreover, $ws[1,i]$ has only $m-1$ elements $\geq j$,
so $y\not\leq ws$. The proof proceeds as
in Case 1, with $i+1$ in place of $i-1$.
\newline
\textbf{Case 3: }$s\neq s_{i-1},s_i$. The
fundamental recursive property of Kazhdan-Lusztig polynomials 
(\cite[Section 7.11, (23)]{humphreys})
tells us that
\beq \label{kleqn}
P_{y,w}=P_{y',ws}+qP_{y's,ws}-
\sum_{\substack{y\leq z<ws\\zs<z}} \mu(z,ws)\,
q^{(\ell(w)-\ell(z))/2}P_{y,z},
\eeq
where $\mu(z,ws)$
is the coefficient of $q^{(\ell(ws)-\ell(z)-1)/2}$ in $P_{z,ws}$, and
$y'$ is the minimum of $y$ and $ys$ in Bruhat order. 
All the nonzero Kazhdan-Lusztig polynomials
involved in the right-hand side are indexed by elements of 
the interval $[y',ws]$,
for which $i$ is cancellable. By the induction hypothesis,
they can all be replaced by the analogous polynomials for the interval
$[(y')^{\hat{i}},(ws)^{\hat{i}}]$, and the result follows.
\epf

We now recall (and
extend slightly) the matrix notation used in \cite{tworemarks}. 
Let $(b_i)_{i\in [1,n]}$ be
an $n$-tuple of nonnegative integers whose sum is $d$, and let
$(c_j)_{j\in[1,n']}$ be an $n'$-tuple of nonnegative integers
whose sum is also $d$.
To avoid notational clutter, we make the convention for the rest of 
this section that the range of the 
variables $i$ and $i'$ will be $[1,n]$ unless otherwise specified, and that
of the variables $j$ and $j'$ will be $[1,n']$. We will use boldface letters
such as $\m$ and $\m'$ for the $(n\times n')$-matrices whose entries are
written with the corresponding ordinary letters $m_{i,j}$ and $m_{i,j}'$.
Let $M_{(b_i);(c_j)}$ be 
the set of all $(n\times n')$-matrices $\m$ satisfying:
\ben
\item $m_{i,j}\in\N$, for all $i,j$,
\item $\sum_{j}m_{i,j}=b_i$, for all $i$,\text{ and } 
\item $\sum_{i}m_{i,j}=c_j$, for all $j$.
\een
If any $b_i$ or $c_j$ is $0$, the corresponding row or column 
must always be zero
and is therefore irrelevant, but it will be convenient
to allow this possibility. We will use an obvious notation for the sums
of various sectors of a matrix:
\[ m_{\leq i,\geq j}=\sum_{\substack{i'\leq i\\j'\geq j}}m_{i',j'},\
m_{\leq i,j}=\sum_{i'\leq i}m_{i',j},\ 
m_{i,\geq j}=\sum_{j'\geq j}m_{i,j'}, \]
and similarly $m_{\geq i,\leq j}$, etc. Note that for $\m\in M_{(b_i);(c_j)}$,
\beq \label{rewriteeqn}
\begin{split} 
m_{\geq i,\leq j}&=c_1+c_2+\cdots+c_j-m_{\leq i-1,\leq j}\\
&=c_1+\cdots+c_j-b_1-\cdots-b_{i-1}+m_{\leq i-1,\geq j+1}.
\end{split}
\eeq

The matrices in $M_{(b_i);(c_j)}$ 
parametrize double cosets of $S_d$ with respect to 
certain parabolic subgroups. 
Namely, write $[1,d]$ as the disjoint union of blocks
$B_1,\cdots,B_n$ such that all elements of 
$B_i$ are less than all elements of $B_{i+1}$,
and $|B_i|=b_i$. (Because we are allowing some $b_i$ to be zero, some
of these blocks could be empty.)
Let $S_{(b_i)}$ be the subgroup
of $S_d$ which preserves each $B_i$ separately; this is a parabolic
subgroup isomorphic to
$S_{b_1}\times\cdots\times S_{b_n}$. Similarly define blocks
$C_j$ of sizes $c_j$, and the parabolic subgroup $S_{(c_j)}$.
We define a surjective map $\psi:S_d\to M_{(b_i);(c_j)}$ by
\[ \psi(w)_{i,j}=|w(B_i)\cap C_j|. \]
The fibres of $\psi$ are exactly the
double cosets $S_{(c_j)}wS_{(b_i)}$, so $\psi$ induces a bijection
$S_{(c_j)}\leftdiv S_d\,/\, S_{(b_i)}\leftrightarrow M_{(b_i);(c_j)}$.
For $\m\in M_{(b_i);(c_j)}$, let $w_{\m}\in S_d$ be the 
longest element in
the corresponding double coset. 

Note that in the case when $n=n'=d$ and all $b_i=c_j=1$, the parabolic
subgroups are trivial, and we have merely passed from elements of $S_d$
to the corresponding permutation matrices (or their transposes, depending
on your convention).
In general,
the permutation $w_\m$ can be constructed from the matrix $\m$
as follows: assuming that the images of $B_{i'}$ for $i'<i$ have
been determined, we send successive various-sized sub-blocks of $B_i$
to the various $C_j$s, according to the entries of the 
$i$th row of $\m$ read from right to left. Within each sub-block,
we successively take the largest element of $C_j$ still unused.
More formally, if $a$ is the $s$th element of $B_i$,
then $w_\m(a)\in C_j$ where $j$ is maximal such that $m_{i,\geq j}
\geq s$. Specifically, $w_\m(a)$ is the $t$th largest element of
$C_j$ where
\beq \label{teqn}
t=m_{\leq i-1,j} + s - m_{i,\geq j+1}.
\eeq
\bex
Take $d=9$, $n=n'=4$, and define $b_i$, $c_j$ so that
\bes
\begin{split} 
&B_1=\{1\},\ B_2=\{2,3,4,5\},\ B_3=\{6,7,8\},\ B_4=\{9\},\\
&C_1=\{1,2\},\ C_2=\{3,4,5\},\ C_3=\{6,7,8\},\ C_4=\{9\}.
\end{split}
\ees
Let us construct $w_\m$ where
\[ \m=\begin{pmatrix}1&0&0&0\\1&2&1&0\\0&1&2&0\\0&0&0&1\end{pmatrix}. \]
The first row tells us that $w_\m(1)$ is an element of $C_1$; we take the
largest element, namely $2$. The second row tells us that $w_\m(B_2)$ consists
of one element of $C_3$, two elements of $C_2$, and one of $C_1$, in that
order. Taking the largest elements not yet used, we set $w_\m(2)=8$,
$w_\m(3)=5$, $w_\m(4)=4$, and $w_\m(5)=1$. Continuing in this way, we see
that $w_\m$ is the permutation $285417639$ (in `one-line' notation).
\eex

We define a length function $\ell:M_{(b_i);(c_j)}\to\N$ by
$\ell(\m)=\ell(w_\m)$, and a partial order on
$M_{(b_i);(c_j)}$ by
\[ \m\leq\m'\Leftrightarrow w_\m\leq w_{\m'}. \]
These can be described as follows.
\bpr \label{finmatbruhatprop}
Let $\m,\m'\in M_{(b_i);(c_j)}$.
\ben
\item $\ell(\m)=\sum_{i,j}
m_{i,j}m_{\leq i,\geq j}-\sum_{i,j} \binom{m_{i,j}+1}{2}$.
\item $\m\leq\m'$ if and only if, for all $i,j$,
\[ m_{\leq i,\geq j}\leq m_{\leq i,\geq j}'. \]
\item $\m\leq\m'$ if and only if, for all $i,j$,
\[ m_{\geq i,\leq j}\leq m_{\geq i,\leq j}'. \]
\een
\epr
\bpf
Let $a$ be the largest element of $B_i\cap w_\m^{-1}(C_j)$.
Then for $1\leq k\leq m_{i,j}$, $a-k+1$ is the $k$th largest element
of $B_i\cap w_\m^{-1}(C_j)$. Clearly
\beq \inv_{a-k+1}(w_\m)=m_{\leq i,\geq j} - k.\eeq
Summing this over all $i,j$, 
and $1\leq k\leq m_{i,j}$ gives (1).
To prove (2), fix $i$ and $j$, and let $b$ be the largest element of
$\cup_{i'\leq i}B_{i'}$ and $c$ the smallest element of 
$\cup_{j'\geq j}C_{j'}$. If $\m\leq\m'$, then by Proposition 
\ref{finbruhatprop}, we have
\[ |\{a\leq b\,|\,w_\m(a)\geq c\}|\leq
|\{a\leq b\,|\,w_{\m'}(a)\geq c\}|, \]
which exactly says that $m_{\leq i,\geq j}\leq m_{\leq i,\geq j}'$. 
Conversely, suppose we know that
$m_{\leq i,\geq j}\leq m_{\leq i,\geq j}'$ and
$m_{\leq i-1,\geq j}\leq m_{\leq i-1,\geq j}'$. For all $1\leq k\leq b_i$, 
we have
\beq |\{a\leq b-k+1\,|\,w_\m(a)\geq c\}|=
\max\{m_{\leq i,\geq j}-k+1,
m_{\leq i-1,\geq j}\}
\eeq
and similarly for $\m'$, so our assumption implies
\[ |\{a\leq b-k+1\,|\,w_\m(a)\geq c\}|\leq
|\{a\leq b-k+1\,|\,w_{\m'}(a)\geq c\}|. \]
Combining these statements for all $j$ tells us that for all $m$,
the $m$th largest element of $w_\m[1,b-k+1]$ lies in a block $C_j$
prior or equal to that containing the $m$th largest element
of $w_{\m'}[1,b-k+1]$. Remembering how $w_\m$ and $w_{\m'}$ are
constructed from $\m$ and $\m'$, this implies that for all $m$,
the $m$th largest element of $w_\m[1,b-k+1]$ is less than or equal to
the $m$th largest element of $w_{\m'}[1,b-k+1]$. Letting $i$ and $k$ vary,
we get $w_\m\leq w_{\m'}$ by Proposition \ref{finbruhatprop}, so (2) is proved.
One way to deduce (3) from (2) is to use \eqref{rewriteeqn}.
Another way is to recall that
$w_\m\leq w_{\m'}$ if and only if $w_\m^{-1}\leq w_{\m'}^{-1}$;
clearly the inverse of $w_\m$ is the permutation $w_{\m^t}$
associated to the transpose matrix $\m^t\in M_{(c_j);(b_i)}$, and the condition
in (3) is the transpose of the condition in (2).
\epf

We can also define Kazhdan-Lusztig polynomials indexed by pairs of elements of
$M_{(b_i);(c_j)}$: $P_{\m,\m'}=P_{w_\m,w_{\m'}}$.
By definition we have
\beq \label{uppertriangeqn}
P_{\m,\m'}\neq 0\Rightarrow \m\leq \m',\text{ and } P_{\m,\m}=1.
\eeq
So the inverse matrix 
$(P_{\m,\m'}^{\langle -1\rangle})_{\m,\m'\in M_{(b_i);(c_j)}}$ of
$(P_{\m,\m'})_{\m,\m'\in M_{(b_i);(c_j)}}$ has entries in $\Z[q]$
which also satisfy \eqref{uppertriangeqn}. In fact, we can express
these entries in terms of those of the original matrix, as follows.
Recall the Kazhdan-Lusztig inversion
formula (\cite[Section 7.14, (24)]{humphreys}):
\beq \label{symminversioneqn}
\sum_{x\in S_d} \varepsilon(xy)\, P_{xw_0^{(d)},yw_0^{(d)}}\, P_{x,w}
=\delta_{y,w},
\eeq
where $\varepsilon(z)=(-1)^{\ell(z)}$ and $w_0^{(d)}$ is the longest
element of $S_d$. Using the fact that 
$P_{x,w_{\m''}}=P_{x',w_{\m''}}$ for all $x'\in S_{(c_j)}x S_{(b_i)}$
(\cite[Section 7.14, Corollary]{humphreys}),
we get
\beq \label{inversioneqn}
P_{\m,\m'}^{\langle -1\rangle}=\sum_{x\in S_{(c_j)}w_{\m'}S_{(b_i)}}
\negthickspace
\varepsilon(xw_\m)\, P_{xw_0^{(d)},w_\m w_0^{(d)}}.
\eeq

A general Kazhdan-Lusztig polynomial $P_{y,w}$, $y,w\in S_d$, can be expressed
in the form $P_{\m,\m'}$ in various ways. The most trivial takes $n=n'=d$ and 
all $b_i=c_j=1$, so that there is no difference between permutations and
matrices.
At the other extreme of usefulness, we can take $(B_i)$ to be the 
collection consisting of
the maximal intervals on which $w$ is decreasing, 
and $(C_j)$ the same thing for $w^{-1}$.
With these choices, $w$ is clearly the longest element in its double coset
$S_{(c_j)}w S_{(b_i)}$, so
$P_{y,w}$ depends only on the double coset of $y$; in other words,
$P_{y,w}=P_{\psi(y),\psi(w)}$ where $\psi:S_d\to M_{(b_i);(c_j)}$ is 
as above.
\bex
Let $y=128456379$ and $w=587429316$ in $S_9$. The blocks $B_i$ and $C_j$
determined by $w$ are exactly those used in the previous example. Indeed,
$w=w_{\m'}$ where
\[ \m'=\begin{pmatrix}0&1&0&0\\1&1&2&0\\1&1&0&1\\0&0&1&0\end{pmatrix}. \]
Now $\psi(y)$ is the matrix $\m$ from the previous example, so 
the permutation $w_\m=285417639$ found there is the longest element in the
double coset $S_{(c_j)}yS_{(b_i)}$. Using the criteria in Proposition
\ref{finmatbruhatprop}, it is easy to check that $\m\leq\m'$.
The above principle means in this case that
\[ P_{y,w}=P_{\m,\m'}=P_{w_\m,w}. \]
The advantage of the latter form is that $2$ is cancellable for $[w_\m,w]$.
Since $w_\m^{\hat{2}}=25417638$ and $w^{\hat{2}}=57428316$, we get
$P_{y,w}=P_{25417638,57428316}$.
\eex

In order to be able to perform such a cancellation directly on matrices,
we note the following.
\bpr \label{finmatcancprop}
Let $\m\in M_{(b_i);(c_j)}$, $a\in B_i\cap w_\m^{-1}(C_j)$.
Let $\mathbf{e}$ be the matrix with
$e_{i,j}=1$, all other entries zero.
\ben
\item $w_\m^{\hat{a}}=w_{\m-\mathbf{e}}$.
\item $\ell(\m)-\ell(\m-\mathbf{e})$ equals each of the following:
\bes
\begin{split}
m_{\leq i,\geq j}+m_{\geq i,\leq j}-m_{i,j}-1&=
m_{\leq i-1,\geq j}+m_{\geq i+1,\leq j}+b_i-1\\
&=m_{\leq i,\geq j+1}+m_{\geq i,\leq j-1}+c_j-1.
\end{split}
\ees
\een
\epr
\bpf
Part (1) is clear from the explicit construction of $w_\m$ given above,
and (2) follows easily from (1) of Proposition \ref{finmatbruhatprop}.
\epf
\bdf
If $\m\leq\m'$ in $M_{(b_i);(c_j)}$, we say that $(i,j)\in[1,n]
\times[1,n']$
is \emph{cancellable} for the interval $[\m,\m']$ if
\ben
\item $m_{i,j}\geq 1$.
\item $m_{\leq i-1,\geq j}=m_{\leq i-1,\geq j}'$, or equivalently
$m_{\geq i,\leq j-1}=m_{\geq i,\leq j-1}'$.
\item $m_{\leq i,\geq j+1}=m_{\leq i,\geq j+1}'$, or equivalently
$m_{\geq i+1,\leq j}=m_{\geq i+1,\leq j}'$.
\een
These equivalences follow from \eqref{rewriteeqn}.
\edf
\bpr \label{finreadyprop}
Suppose that $\m\leq\m'$ in $M_{(b_i);(c_j)}$ and $(i,j)$ is 
cancellable for $[\m,\m']$. Let $\mathbf{e}$ be the matrix with
$e_{i,j}=1$, all other entries zero.
\ben
\item For any $\m^1\in[\m,\m']$,
\ben
\item $m_{i,j}^1\geq m_{i,j}$,
\item $m_{\leq i-1,\geq j}^1=m_{\leq i-1,\geq j}$, and
\item $m_{\leq i,\geq j+1}^1=m_{\leq i,\geq j+1}$.
\een
Hence $(i,j)$ is cancellable for any sub-interval of $[\m,\m']$.
\item The map $\m^1\mapsto\m^1-\mathbf{e}$ is an isomorphism of posets
between $[\m,\m']$ and $[\m-\mathbf{e},\m'-\mathbf{e}]$,
which reduces all lengths by the same amount.
\item For any $\m^1,\m^2\in[\m,\m']$, 
$P_{\m^1,\m^2}=P_{\m^1-\mathbf{e},\m^2-\mathbf{e}}$.
\item For any $\m^1,\m^2\in[\m,\m']$, 
$P_{\m^1,\m^2}^{\langle -1\rangle}=
P_{\m^1-\mathbf{e},\m^2-\mathbf{e}}^{\langle -1\rangle}$.
\een
\epr
\bpf
Let $\m^1\in[\m,\m']$. By (2) of Proposition \ref{finmatbruhatprop}, we have
\[ m_{\leq i-1,\geq j}\leq m_{\leq i-1,\geq j}^1\leq m_{\leq i-1,\geq j}'
=m_{\leq i-1,\geq j}, \]
which proves (1b), and (1c) is similar. It follows that
\bes
\begin{split}
m_{i,j}^1-m_{i,j}&=(m_{i,j}^1+m_{\leq i-1,\geq j}^1+m_{\leq i,\geq j+1}^1)
-(m_{i,j}+m_{\leq i-1,\geq j}+m_{\leq i,\geq j+1})\\
&=(m_{\leq i,\geq j}^1+m_{\leq i-1,\geq j+1}^1)-
(m_{\leq i,\geq j}+m_{\leq i-1,\geq j+1})\\
&=(m_{\leq i,\geq j}^1-m_{\leq i,\geq j})+(m_{\leq i-1,\geq j+1}^1-
m_{\leq i-1,\geq j+1})\\
&\geq 0,
\end{split}
\ees
by Proposition \ref{finmatbruhatprop} again. So (1a) is proved. 
Thus $\m^1-\mathbf{e}\in M_{(\widetilde{b_i});
(\widetilde{c_j})}$, where
\[ \widetilde{b_{i'}}=b_{i'}-\delta_{i,i'},\ 
\widetilde{c_{j'}}=c_{j'}-\delta_{j,j'}. \]
Given this, the first part of (2) is obvious from either 
description of the partial order
given in Proposition \ref{finmatbruhatprop}, and the second part from
(2) of Proposition \ref{finmatcancprop}. To prove (3), let
$a$ be the largest element of $B_i\cap w_{\m}^{-1}(C_j)$, i.e.\
the $(m_{i,\geq j})$th element of $B_i$. We want to
show that $a$ is cancellable for $[w_\m,w_{\m'}]$. Let
\[ \delta=m_{i,\geq j}'-m_{i,\geq j}=m_{\leq i,j}'-m_{\leq i,j}
=m_{\leq i,\geq j}'-m_{\leq i,\geq j}\geq 0. \]
(These are equal because $(i,j)$ is cancellable, and nonnegative
because $\m\leq\m'$.) By the above chain of equalities applied to
$\m^1=\m'$, we have 
\[ m_{i,j}'-m_{i,j}=\delta+(m_{\leq i-1,\geq j+1}'-
m_{\leq i-1,\geq j+1})\geq\delta, \]
so
\[ m_{i,\geq j+1}'=m_{i,\geq j}'-m_{i,j}'=m_{i,\geq j}+\delta-m_{i,j}'<
m_{i,\geq j}\leq m_{i,\geq j}', \]
which means that $w_{\m'}(a)\in C_j$. From 
\eqref{teqn} we see that $w_{\m'}(a)=w_{\m}(a)$. Moreover,
\[ \inv_{a}(w_\m)=m_{\leq i,\geq j} - 1
=m_{\leq i,\geq j}' - \delta- 1
=\inv_{a}(w_{\m'}), \]
so $a$ is cancellable for $[w_\m,w_{\m'}]$. In particular, for any
$\m^1\in[\m,\m']$, $w_{\m^1}(a)\in C_j$, which by (1) of Proposition
\ref{finmatcancprop} implies $w_{\m^1}^{\hat{a}}=w_{\m^1-\mathbf{e}}$.
Then part (3) follows from (3) of Proposition 
\ref{fincancprop}, and part (4) follows formally from parts (2) and (3).
\epf
\bex
With matrices $\m$, $\m'$ defined as in previous examples,
$(2,3)$ is cancellable for $[\m,\m']$, corresponding to the fact that
$2$ is cancellable for $[w_\m,w_{\m'}]$. Performing the cancellation
directly on the matrices, we get
\[ \m-\mathbf{e}=
\begin{pmatrix}1&0&0&0\\1&2&0&0\\0&1&2&0\\0&0&0&1\end{pmatrix},\
\m'-\mathbf{e}=
\begin{pmatrix}0&1&0&0\\1&1&1&0\\1&1&0&1\\0&0&1&0\end{pmatrix}. \]
The reader can check that these matrices correspond to the permutations
$25417638$ and $57428316$ found earlier.
\eex
\section{Nilpotent orbits of the linear quiver}
We now return to the set-up of the first part of the introduction,
so $V$ is a $d$-dimensional $\Z$-graded vector space, with
$d_i=\dim V_i$. For convenience, we adjust the grading so that
$d_i\neq 0\Rightarrow i\in [1,n]$, for some positive integer $n$
(so we are effectively considering the linear quiver of type $A_n$).
Throughout this section, the variables $i,j$ range
over $[1,n]$ unless otherwise specified. 

We saw in \S1 that the $G_V$-orbits in $\CalN_V$ are in bijection with the
set $M_{(d_i)}$ of multisegments in which $i$ occurs $d_i$ times
as an element of a segment.
Following \cite{tworemarks}, we change this parametrization by multisegments
to a parametrization by matrices. We identify each $\m\in M_{(d_i)}$ 
with the $(n\times n)$-matrix $(m_{i,j})$, where
\beq \label{conditioneqn}
m_{i,j}=\left\{\begin{array}{cl}
\text{multiplicity of the segment $[i,j]$},&\text{ if $i\leq j$,}\\
\text{number of segments $[k,l]$ where $k\leq j$, $l\geq i$},
&\text{ if $j=i-1$,}\\
0,&\text{ if $j<i-1$.}
\end{array}\right.
\eeq
It is clear that this matrix lies in the set
$M_{(d_i);(d_j)}$, as defined in the previous section. So we have identified
$M_{(d_i)}$ with a subset of $M_{(d_i);(d_j)}$, which can be described 
as follows.
\bpr \label{finidealprop}
Let $M_{(d_i)}'=\{\m\in M_{(d_i);(d_j)}\,|\,m_{i,j}=0,\,\forall j<i-1\}$.
\ben
\item $M_{(d_i)}'$ is a lower ideal of the poset $M_{(d_i);(d_j)}$.
\item If $\m\in M_{(d_i)}'$, then for all $i\geq j$,
$m_{\leq i,\geq j}=d_j+d_{j+1}+\cdots+ d_i$.
\item For $\m,\m'\in M_{(d_i)}'$, $\m\leq\m'$ if and only if
$m_{\leq i,\geq j}\leq m_{\leq i,\geq j}'$ for all $i<j$.
\item If $\m\in M_{(d_i)}'$, then for all $i\in [2,n]$,
$m_{i,i-1}=m_{\leq i-1,\geq i}$.
\item  $M_{(d_i)}'=M_{(d_i)}$.
\een
\epr
\bpf
An element $\m\in M_{(d_i);(d_j)}$ lies in $M_{(d_i)}'$ if and only
if $m_{\geq i,\leq j}=0$ for all $i,j$ such that $j<i-1$, so (1) follows
from (3) of Proposition \ref{finmatbruhatprop}.
For (2), since $m_{\geq i+1,\leq j-1}=0$, \eqref{rewriteeqn} gives
\[ m_{\leq i,\geq j}=d_1+\cdots+d_i-d_{1}-\cdots -d_{j-1}=d_j+\cdots+d_i, \]
as required. Part (3) then follows from (2) of Proposition 
\ref{finmatbruhatprop}.
For (4), we have
\[ m_{i,i-1}=d_i-m_{i,\geq i}=m_{\leq i,\geq i}-m_{i,\geq i}
=m_{\leq i-1,\geq i}. \]
From (4) and the $i=j$ case of (2) 
it follows that every matrix in $M_{(d_i)}'$ arises
from a multisegment in $M_{(d_i)}$ by the rule \eqref{conditioneqn}, 
whence (5).
\epf
As mentioned in the introduction, the identification of $M_{(d_i)}$
with $M_{(d_i)}'$ is a poset isomorphism: the geometrically-defined
partial order $\preceq$ on $M_{(d_i)}$ is the restriction of the
partial order $\leq$ on $M_{(d_i);(d_j)}$. This is part of
Zelevinsky's result
\cite[Corollary 1]{tworemarks}, which we can state 
(with some supplementary detail) as follows.
\bth \label{zelevinskythm}
Let $\m,\m'\in M_{(d_i)}$.
\ben
\item $\dim\CalO_\m=\ell(\m)-\sum_{i}\binom{d_i}{2}$.
\item $\CalO_\m\subseteq
\overline{\CalO_{\m'}}\Leftrightarrow \m\leq \m'$.
\item $\CalH^i IC(\overline{\CalO_{\m'}})=0$ for $i$ odd. 
\item $IC_{\m,\m'}=P_{\m,\m'}$.
\item $IC_{\m,\m'}^{\langle -1\rangle}=P_{\m,\m'}^{\langle -1\rangle}$.
\een
\eth
\bpf
For reference in \S5, we recall Zelevinsky's proof.
Define the partial flag variety $\B_{(d_i)}$ to be the set of collections of 
subspaces $(W_i)_{i\in[0,n]}$ of $V$
such that $W_0=0$, and for all $i\in [1,n]$,
$W_{i-1}\subset W_{i}$ and $\dim W_{i}/W_{i-1}=d_i$; this is naturally
a nonsingular projective variety of dimension 
$\binom{d}{2}-\sum_i\binom{d_i}{2}$.
We define a `base-point' $(U_i)$ in $\B_{(d_i)}$ by
$U_i=V_1\oplus\cdots\oplus V_i$.
Relative to this base-point, $\B_{(d_i)}$ decomposes into
Schubert cells $\B_\m$ for $\m\in M_{(d_i);(d_j)}$.
Explicitly,
$\B_\m$ consists of those $(W_i)$ such that for all $i,j\in[1,n]$,
\[ \dim \frac{W_i\cap U_j}{W_i\cap U_{j-1}+
W_{i-1}\cap U_j} =m_{i,j}. \]
The analogues of (1)--(4) for these Schubert cells
(for all of $M_{(d_i);(d_j)}$) are well known.
Let $\B_{(d_i)}'$ be the closed subvariety of $\B_{(d_i)}$
defined by requiring
$W_i\supset U_{i-1}$; from the description of $M_{(d_i)}$ as $M_{(d_i)}'$,
it is easy to see that
$\B_{(d_i)}'=\bigcup_{\m\in M_{(d_i)}}\B_\m$.

Now we define a morphism $\CalN_V\to \B_{(d_i)}':
\varphi\mapsto(W_i(\varphi))$ by the rule
\[ W_i(\varphi)=U_{i-1}\oplus\{v+\varphi(v)+\varphi^2(v)+\cdots
+\varphi^{n-i}(v)\,|\,v\in V_i\}. \]
An easy check shows that this morphism maps
$\CalO_\m$ into $\B_\m$ for all $\m\in M_{(d_i)}$.
Moreover, it gives an isomorphism between
$\CalN_V$ and the open subvariety of $\B_{(d_i)}'$
defined by requiring
\[ W_i\cap\bigoplus_{i'>i}V_{i'}=0,\ \forall i\in[1,n]. \]
Hence each $\CalO_\m$ is embedded as an open subvariety of the Schubert cell
$\B_\m$, and (1)--(4) follow. Since $M_{(d_i)}$ is a lower ideal
of $M_{(d_i);(d_j)}$, (5) is an automatic consequence of (4).
\epf

In view of (3) of Proposition \ref{finidealprop}, part (2) of Theorem
\ref{zelevinskythm} says that $\CalO_\m\subseteq
\overline{\CalO_{\m'}}$ if and only if for all $i<j$,
$m_{\leq i,\geq j}\leq m_{\leq i,\geq j}'$. Now if $\varphi\in\CalO_\m$, then
for $i\leq j$, $m_{\leq i,\geq j}=\rk\varphi^{j-i}|_{V_{i}}$. So we recover the
well-known fact that $\CalO_\m\subseteq
\overline{\CalO_{\m'}}$ if and only if for all $i<j$,
$\rk\varphi^{j-i}|_{V_i}\leq\rk(\varphi')^{j-i}|_{V_i}$ 
for any $\varphi\in\CalO_{\m}$
and $\varphi'\in\CalO_{\m'}$. (Of course the ``only if'' direction is
obvious.)
We can define an element $\mmax\in M_{(d_i)}$ uniquely by the
requirement that for $i\leq j$, $m_{\leq i,\geq j}^{\mathrm{max}}$ equals
the maximum possible rank, namely $\min\{d_i,d_{i+1},\cdots,d_j\}$.
It follows that $\m\leq\mmax$ for all $\m\in M_{(d_i)}$, and
the orbit $\CalO_{\mmax}$ is dense in $\CalN_V$.

As foreshadowed in the introduction, Theorem \ref{zelevinskythm} is
only one of many possible ways to express a particular $IC_{\m,\m'}$
as a Kazhdan-Lusztig polynomial: the below-diagonal entries prescribed
by \eqref{conditioneqn} correspond to one particular choice of
``empty segments''. A more general statement is the following.
\bth \label{finmainthm}
Let $b_1,\cdots,b_n,c_1,\cdots,c_n\in\N$ be such that
\[ b_1=d_1,\ c_n=d_n,\text{ and }d_i-b_i=d_{i-1}-c_{i-1},\,
\forall i\in [2,n]. \]
Define an $(n\times n)$-matrix $\mathbf{a}$ by
\[ a_{i,j}=\left\{\begin{array}{cl}
d_i-b_i,&\text{ if $j=i-1$,}\\
0,&\text{ otherwise.}
\end{array}\right. \]
Let $M_{(d_i)}^{(b_i);(c_j)}=\{\m\in M_{(d_i)}\,|\,m_{i,i-1}\geq d_i-b_i,\,
\forall i\in [2,n]\}$.
\ben
\item $M_{(d_i)}^{(b_i);(c_j)}$ is an upper ideal of $M_{(d_i)}$.
\item The map $\m\mapsto\m-\mathbf{a}$ is an isomorphism of posets between
$M_{(d_i)}^{(b_i);(c_j)}$ and $\{\widetilde{\m}\in M_{(b_i);(c_j)}\,|\,
\widetilde{m}_{i,j}=0,\,\forall j< i-1\}$.
\item For any $\m',\m''\in M_{(d_i)}^{(b_i);(c_j)}$, 
$IC_{\m',\m''}=P_{\m'-\mathbf{a},\m''-\mathbf{a}}$.
\item For any $\m',\m''\in M_{(d_i)}^{(b_i);(c_j)}$, 
$IC_{\m',\m''}^{\langle -1\rangle}=
P_{\m'-\mathbf{a},\m''-\mathbf{a}}^{\langle -1\rangle}$.
\een
\eth
\bpf
For $\m\in M_{(d_i)}$, $m_{i,i-1}$ equals $m_{\geq i,\leq i-1}$, so (1) follows
from (3) of Proposition \ref{finmatbruhatprop}. For (2), the fact that the
given map is a bijection is obvious, and by Proposition \ref{finmatbruhatprop}
it preserves partial orders. Now this map is the composition of
maps of the form $\m\mapsto\m-\mathbf{e}$ as in Proposition \ref{finreadyprop}
and their inverses $\m\mapsto\m+\mathbf{e}$, where the positions which
are being altered are all of the form $(i,i-1)$. Since all matrices involved
have zero entries in positions $(i,j)$ where $j<i-1$, conditions (2) and (3)
of the definition of cancellability always hold. So Proposition 
\ref{finreadyprop} implies (3) and (4) with $P_{\m',\m''}$ and
$P_{\m',\m''}^{\langle -1\rangle}$ in place of
$IC_{\m',\m''}$ and $IC_{\m',\m''}^{\langle -1\rangle}$, and
Theorem \ref{zelevinskythm} gives the result.
\epf
\bcr \label{finmaincor}
For $\m\in M_{(d_i)}$, define an $(n\times n)$-matrix $\m^-$ by
\[ m_{i,j}^-=\left\{\begin{array}{cl}
m_{i,j},&\text{ if $j=i-1$,}\\
0,&\text{ otherwise.}
\end{array}\right. \]
Let $\langle\m\rangle=\{\m'\in M_{(d_i)}\,|\,\m\leq\m'\}=[\m,\mmax]$.
\ben
\item For all $\m'\in\langle\m\rangle$, 
$\m'-\m^-\in M_{(d_i-m_{i,i-1});(d_j-m_{j+1,j})}$.
\item The map $\m'\mapsto\m'-\m^-$ is an isomorphism of posets between
$\langle\m\rangle$ and $[\m-\m^-,\mmax-\m^-]$.
\item For any $\m',\m''\in\langle\m\rangle$, 
$IC_{\m',\m''}=P_{\m'-\m^-,\m''-\m^-}$.
\item For any $\m',\m''\in\langle\m\rangle$, 
$IC_{\m',\m''}^{\langle -1\rangle}=P_{\m'-\m^-,\m''-\m^-}^{\langle -1\rangle}$.
\een
\ecr
\bpf
Apply Theorem \ref{finmainthm} with
\[ b_i=d_i-m_{i,i-1},\ c_{i-1}=d_{i-1}-m_{i,i-1} \]
for all $i\in [2,n]$, and restrict to the upper ideal $\langle\m\rangle$
of $M_{(d_i)}^{(b_i);(c_j)}$.
\epf
Note that the polynomials $P_{\m'-\m^-,\m''-\m^-}$ in (3) are
Kazhdan-Lusztig polynomials for $S_{k(\m)}$,
where $k(\m)$ is the number of segments of $\m$, which is also the sum of the
entries of $\m-\m^-$.

Finally, we have to connect Theorem \ref{finmainthm}
to the notation of the introduction,
in order to prove Theorem \ref{introfinmainthm}. We have elements
$\lambda,\mu\in D_k$; we can clearly assume that all $\lambda_s-s,\mu_s-s+1$
for $1\leq s\leq k$ lie in $[1,n]$. Define
\[ b_i=|\{s\,|\,\mu_s-s+1=i\}|,\ c_j=|\{s\,|\,\lambda_s-s=j\}|. \]
Then the subgroups $S_{(b_i)}$ and $S_{(c_j)}$ of $S_k$
are exactly the conjugates
under $w_0^{(k)}$ of the dot stabilizers $W_\mu$ and $W_\lambda$
(this reversal comes about because the sequences $(\lambda_s-s)$ and
$(\mu_s-s)$ are decreasing).
The map $\psi:S_k\to M_{(b_i);(c_j)}$
as defined in the previous section satisfies
\beq \label{finpsieqn}
\psi(w_0^{(k)}ww_0^{(k)})_{i,j}=|\{s\,|\,\mu_s-s+1=i,
\lambda_{w(s)}-w(s)=j\}|.
\eeq
So $S_k[\lambda,\mu]=\{w\in S_k\,|\,\psi(w_0^{(k)}ww_0^{(k)})_{i,j}=0,\,\forall
j<i-1\}$, which shows that it is indeed a lower ideal of $S_k$.
Moreover, $w\mapsto \psi(w_0^{(k)}ww_0^{(k)})$ gives an 
isomorphism of posets between
$S_k[\lambda,\mu]^\circ$ and
$\{\widetilde{\m}\in M_{(b_i);(c_j)}\,|\, \widetilde{m}_{i,j}=0,\,\forall
j<i-1\}$, and the polynomials attached to these posets correspond,
since
\beq 
P_{w,w'}=P_{w_0^{(k)}ww_0^{(k)},w_0^{(k)}w'w_0^{(k)}}
=P_{\psi(w_0^{(k)}ww_0^{(k)}),\psi(w_0^{(k)}w'w_0^{(k)})}
\eeq
for all $w,w'\in S_k[\lambda,\mu]^\circ$. Now as in Theorem 
\ref{introfinmainthm}, assume that $\lambda/\mu\in M_{(d_i)}$; it follows
immediately that $(b_i)$ and $(c_j)$ satisfy the conditions of
Theorem \ref{finmainthm}. By \eqref{finpsieqn},
for all $w\in S_k[\lambda,\mu]^\circ$,
the multisegment $\lambda/(w\cdot\mu)$ when viewed as a matrix has
the same diagonal and above-diagonal entries as $\psi(w_0^{(k)}ww_0^{(k)})$;
hence 
$\lambda/(w\cdot\mu)=\psi(w_0^{(k)}ww_0^{(k)})+\mathbf{a}$
where $\mathbf{a}$ is as in Theorem \ref{finmainthm}.
Thus Theorem \ref{introfinmainthm} follows from Theorem \ref{finmainthm}.
\section{Cancellation for the affine symmetric group}
We now want to extend the results of \S2 to the affine symmetric group.
Again fix a positive integer $d$.
Let $\Sdh$ be the group of permutations $w$ of the set $\Z$ such that
$w(i+d)=w(i)+d$, for all $i\in\Z$.
An element $w\in\Sdh$ is determined by its \emph{window}
$(w(1),w(2),\cdots,w(d))$,
which can be any collection of representatives of the congruence classes
mod $d$, in any order. The subgroup of $\Sdh$ which preserves $[1,d]$
is clearly isomorphic to $S_d$.

The group $\Sdh$ is the `extended' affine symmetric group: it 
can be written as a semi-direct
product $\langle\tau\rangle\ltimes\Sdt$, where
\[ \Sdt=\{w\in\Sdh\,|\,\sum_{i=1}^{d} w(i)=\sum_{i=1}^{d} i\} \]
is the actual affine symmetric group, and $\tau$ is the element of infinite 
order sending $i$ to $i+1$ for all $i\in\Z$. In general,
$w\in\tau^{a(w)}\Sdt$ where
\[ a(w)=\frac{1}{d}(\sum_{i=1}^{d} w(i)-\sum_{i=1}^{d} i). \]
Note that for any $i\in\Z$, the set
$w(-\infty,i]$ can be obtained from $(-\infty,i+a(w)]$ by changing
finitely many elements (keeping distinctness). In other words, for $m$
sufficiently large, the $m$th largest element in $w(-\infty,i]$ is
$i+a(w)-m+1$.

If $d=1$, $\Sdt$ is the trivial group and $\Sdh=\langle\tau\rangle$.
If $d\geq 2$, we define
$s_i\in\Sdt$ for all $i\in\Z$ by
\[ s_i(j)=\left\{\begin{array}{cl}
j+1,&\text{ if $j\equiv i$ mod $d$,}\\
j-1,&\text{ if $j\equiv i+1$ mod $d$,}\\
j,&\text{ otherwise.}
\end{array}\right. \]
Thus $s_i=s_{i'}$ iff $i\equiv i'$ mod $d$. It is well known that
$s_0, s_1,\cdots,s_{d-1}$ form a set of Coxeter generators for $\Sdt$
of type $\widetilde{A_{d-1}}$. They thus determine a length function
$\ell:\Sdt\to\N$, a Bruhat order $\leq$ on $\Sdt$, and Kazhdan-Lusztig
polynomials $P_{y,w}\in\N[q]$ for $y,w\in\Sdt$, all of which are
invariant under conjugation by $\tau$. We extend these to $\Sdh$
in the standard way:
\bes
\begin{split}
\ell(w)&=\ell(\tau^{-a(w)}w),\\
y\leq w&\Leftrightarrow a(y)=a(w), \tau^{-a(y)}y\leq\tau^{-a(w)}w,\\
P_{y,w}&=\left\{\begin{array}{cl}
P_{\tau^{-a(y)}y,\tau^{-a(w)}w},&\text{ if $a(y)=a(w)$,}\\
0,&\text{ otherwise.}
\end{array}\right.
\end{split}
\ees

We define inversion statistics as in the finite case:
\[ \inv_{i}(w)=|\{i'<i\,|\,w(i')>w(i)\}|,\
\Inv_i(w)=|\{i'>i\,|\,w(i')<w(i)\}|, \]
for any $w\in\Sdh$ and $i\in\Z$ (these sets are finite,
even though $i'$ runs over $\Z$.)
Clearly $\inv_{i+d}(w)=\inv_{i}(w)$, $\Inv_{i+d}(w)=\Inv_i(w)$, and
\bes
\begin{split}
&\Inv_i(w)-\inv_i(w)\\
&=|(-\infty,w(i)]\setminus w(-\infty,i]|
-|w(-\infty,i]\setminus(-\infty,w(i)]|\\
&=|(-\infty,w(i)]\setminus (-\infty,i+a(w)]|
-|(-\infty,i+a(w)]\setminus(-\infty,w(i)]|\\
&=w(i)-i-a(w).
\end{split}
\ees
The formula for $\ell$ on $\Sdh$ is analogous to that for $S_d$
(see \cite[Proposition 4.1(ii)]{bjornerbrenti}):
\beq \label{lengtheqn}
\ell(w)=\sum_{i\in [1,d]}\inv_i(w)=\sum_{i\in [1,d]}\Inv_i(w).
\eeq
We also have
\beq
ws_i<w\text{ if and only if }w(i)>w(i+1),
\eeq 
and a general description of Bruhat order along the lines of Proposition
\ref{finbruhatprop} (this is a rephrasing of
\cite[Theorem 6.5]{bjornerbrenti}, trivially extended from $\widetilde{S_d}$
to $\widehat{S_d}$):
\bpr \label{bruhatprop}
If $y,w\in\Sdh$,
$y\leq w$ if and only if for all $i\in\Z$,
\[ |\{i'\leq i\,|\,y(i')\geq j\}|\leq |\{i'\leq i\,|\,w(i')\geq j\}|,\
\forall j,\text{ with equality for $j\ll 0$.} \]
\epr
\noindent
In other words, for all positive integers $m$, 
the $m$th largest element in $y(-\infty,i]$
is less than or equal to the $m$th largest element in $w(-\infty,i]$,
with equality for $m\gg 0$. It suffices to check this for $i\in [1,d]$.

The definition of cancellability is identical to the finite case:
\bdf
If $y\leq w$ in $\Sdh$, we say that $i\in\Z$ is \emph{cancellable} for
the interval $[y,w]$ if $y(i)=w(i)$, $\inv_{i}(y)=\inv_{i}(w)$,
and $\Inv_{i}(y)=\Inv_{i}(w)$. (Clearly any two of these conditions imply
the third, and $i$ is cancellable for $[y,w]$ iff $i+d$ is.)
\edf
However, the process of cancellation is not as uniquely defined
as in the finite case: we need to choose order-preserving
bijections $\sigma_{\bar i}:\Z\setminus\bar{i}\to\Z$
for all congruence classes $\bar i$ mod $d$.
Then for any $w\in\Sdh$, we define $w^{\hat{i}}\in\widehat{S_{d-1}}$ by
\[ w^{\hat{i}}=\sigma_{\overline{w(i)}}\circ w\circ\sigma_{\bar{i}}^{-1}. \]
Note that using different $\sigma$'s would have the effect of multiplying
$w^{\hat{i}}$ on left and right by powers of $\tau$. Independently
of the choice, we have
\beq
\ell(w^{\hat{i}})=\ell(w)-\inv_i(w)-\Inv_i(w).
\eeq
\bex
Take $d=3$, $y=\tau s_1 s_2$ and $w=\tau s_2 s_1 s_0 s_2$.
Then $y$ has window $(3,4,2)$ and $w$ has window $(0,7,2)$.
Since $\inv_3(y)=\inv_3(w)=2$, $3$ is cancellable for $[y,w]$.
If we normalize $\sigma_{\bar 2}$ and $\sigma_{\bar 3}$ by requiring
that they preserve $1$, then $y^{\hat 3}$ and $w^{\hat 3}$ are the
elements of $\widehat{S_2}$ with windows $(2,3)$ and $(0,5)$, namely
$\tau$ and $\tau s_1 s_0$.
\eex
We can now extend Proposition \ref{fincancprop} to the affine case.
\bpr \label{cancprop}
Suppose that $i$ is cancellable for $[y,w]$.
\ben
\item For any $x\in [y,w]$, $x(i)=y(i)$ and $\inv_{i}(x)=\inv_{i}(y)$.
Hence $i$ is cancellable for any sub-interval of $[y,w]$.
\item $x\mapsto x^{\hat{i}}$ is an isomorphism of posets between $[y,w]$
and $[y^{\hat{i}},w^{\hat{i}}]$, which reduces all lengths by the same
amount.
\item For any $u,v\in[y,w]$, $P_{u,v}=P_{u^{\hat{i}},v^{\hat{i}}}$.
\een
\epr
\bpf
The proof of part (1) is identical to that of part (1) of Proposition
\ref{fincancprop}, with $[1,i]$ replaced by $(-\infty,i]$ and so on,
and of course using Proposition \ref{bruhatprop} instead of
Proposition \ref{finbruhatprop}. Similarly with part (2), where the inverse map
$[y^{\hat{i}},w^{\hat{i}}]\to[y,w]:x\mapsto\tilde{x}$ is now defined by
\[ \tilde{x}(i')=\left\{\begin{array}{cl}
y(i)+kd,&\text{ if $i'=i+kd$,}\\
\sigma_{\overline{y(i)}}^{-1}(x(\sigma_{\bar i}(i'))),
&\text{ if $i'\not\in\bar i$.}
\end{array}\right. \]
The proof of (3) is also mostly unchanged. Apart from replacing
$[1,i]$ by $(-\infty,i]$ and so on, the only change is that in Case 1,
we need not have $(ys)^{\widehat{i-1}}=y^{\hat{i}}$ and
$(ws)^{\widehat{i-1}}=w^{\hat{i}}$, but rather have
\[ (ys)^{\widehat{i-1}}=\tau^a y^{\hat{i}}\tau^b,\
(ws)^{\widehat{i-1}}=\tau^a w^{\hat{i}}\tau^b\text{ for some
$a,b\in\Z$}, \]
which still implies $P_{(ys)^{\widehat{i-1}},(ws)^{\widehat{i-1}}}=
P_{y^{\hat{i}},w^{\hat{i}}}$ as required.
\epf

We now introduce some affine matrix notation very similar to that
in \cite{aperiodicity}.
Let $(b_i)_{i\in\Z}$ be
a $\Z$-tuple of nonnegative integers, periodic with period $n\geq 1$,
such that $\sum_{i=1}^n b_i=d$; and let $(c_j)_{j\in\Z}$ be another such,
with period $n'\geq 1$, such that $\sum_{j=1}^{n'}c_j=d$. Our notational
convention now is that the range of the variables $i,i',j,j'$ is all
of $\Z$ unless otherwise specified.
Let $M_{(b_i),n;(c_j),n'}$ be 
the set of all $(\Z\times\Z)$-matrices $\m$ satisfying:
\ben
\item $m_{i,j}\in\N$, for all $i,j$,
\item $m_{i+n,j+n'}=m_{i,j}$, for all $i,j$,
\item $\sum_{j}m_{i,j}=b_i$, for all $i$,\text{ and } 
\item $\sum_{i}m_{i,j}=c_j$, for all $j$.
\een
It is easy to see that for $\m\in M_{(b_i),n;(c_j),n'}$, $m_{i,j}=0$ for
$|j-i|\gg 0$; so sums of the form
$m_{i,\geq j}$, $m_{\leq i,j}$, $m_{\leq i,\geq j}$ and 
$m_{\geq i,\leq j}$ are finite. We have the following substitute
for \eqref{rewriteeqn}. For fixed $i$, $m_{\geq i,\leq j_0}=0$
for all $j_0$ sufficiently negative, and for $j$ greater than such $j_0$,
\beq \label{affrewriteeqn}
m_{\geq i,\leq j}=c_{j_0+1}+\cdots+c_j-m_{\leq i-1,\geq j_0+1}
+m_{\leq i-1,\geq j+1}.
\eeq

The matrices in $M_{(b_i),n;(c_j),n'}$
parametrize double cosets of $\Sdh$ with respect to 
proper parabolic subgroups of $\Sdt$. 
Namely, write $\Z$ as the disjoint union of (possibly empty) blocks
$B_i$ such that all elements of $B_i$ are less than all elements of $B_{i+1}$,
and $|B_i|=b_i$.
It follows that $B_{i+n}=B_i+d$. Note that the collection $(B_i)$
is determined by $(b_i)$ up to translation (i.e.\ a power of $\tau$).
Let $S_{(B_i)}$ be the subgroup
of $\Sdh$ which preserves each $B_i$ separately; this is a parabolic
subgroup of $\Sdt$ isomorphic to
$S_{b_1}\times\cdots\times S_{b_n}$. (It is determined by $(b_i)$
up to conjugation by a power of $\tau$.) Similarly define blocks
$C_j$ of sizes $c_j$, and the parabolic subgroup $S_{(C_j)}$.
We define a surjective map $\psi:\Sdh\to M_{(b_i),n;(c_j),n'}$ by
\[ \psi(w)_{i,j}=|w(B_i)\cap C_j|. \]
The fibres of $\psi$ are exactly the
double cosets $S_{(C_j)}wS_{(B_i)}$, so $\psi$ induces a bijection
$S_{(C_j)}\leftdiv\Sdh\,/\, S_{(B_i)}\leftrightarrow M_{(b_i),n;(c_j),n'}$.
For $\m\in M_{(b_i),n;(c_j),n'}$, let $w_{\m}\in\Sdh$ be the 
longest element in
the corresponding double coset. 

The permutation $w_\m$ can be read off the matrix $\m$ by exactly the same
prescription as in the finite case (remembering that $i,i',j,j'$ now
range over all of $\Z$).
\bex
Let $d=7$, $n=2$, $n'=3$, and define
\[ b_i=\left\{\begin{array}{cl}
3,&\text{ if $i\equiv 1$ mod $2$,}\\
4,&\text{ if $i\equiv 0$ mod $2$,}
\end{array}\right.
\quad
c_j=\left\{\begin{array}{cl}
3,&\text{ if $j\equiv 1$ mod $3$,}\\
2,&\text{ if $j\equiv 2$ mod $3$,}\\
2,&\text{ if $j\equiv 0$ mod $3$.}
\end{array}\right. \]
Let $\m\in M_{(b_i),2;(c_j),3}$ be the following matrix:
\[
\begin{array}{cccccccccccccc}
\ddots&&\vdots&&\vdots&&\vdots&&\vdots&&\vdots&&\vdots&\\
&2&0&0&0&0&1&0&0&0&0&0&0&\\
\cdots&0&2&0&0&0&1&0&0&0&0&0&0&\cdots\\
&1&0&0&2&0&0&0&0&1&0&0&0&\\
\cdots&0&0&0&\mathbf{0}&2&0&0&0&1&0&0&0&\cdots\\
&0&0&0&1&0&0&2&0&0&0&0&1&\\
\cdots&0&0&0&0&0&0&0&2&0&0&0&1&\cdots\\
&0&0&0&0&0&0&1&0&0&2&0&0&\\
&&\vdots&&\vdots&&\vdots&&\vdots&&\vdots&&\vdots&\ddots
\end{array}
\]
where the $\mathbf{0}$ is the $(1,1)$ entry.
We choose $(B_i)$ and $(C_j)$ so that $B_1=C_1=\{1,2,3\}$.
The row containing $\mathbf{0}$ tells us that
$w_\m(B_1)$ consists of one element of $C_6=\{13,14\}$ and two of 
$C_2=\{4,5\}$, in that order. Since the $(0,6)$ entry is $1$,
the largest element of 
$C_6$ `has already been used' in $w_\m(B_0)$, so we set $w_\m(1)=13$,
$w_\m(2)=5$, $w_\m(3)=4$. Treating the next row similarly, we find
that $w_\m$ is the element of $\widehat{S_7}$ with window
$(13,5,4,21,10,9,1)$.
\eex

We define a length function $\ell:M_{(b_i),n;(c_j),n'}\to\N$ by
$\ell(\m)=\ell(w_\m)$, and a partial order on
$M_{(b_i),n;(c_j),n'}$ by
\[ \m\leq\m'\Leftrightarrow w_\m\leq w_{\m'}. \]
Since the map $\m\mapsto w_\m$
depends on the choice of $(B_i)$ and $(C_j)$ only modulo 
left and right multiplication
by fixed powers of $\tau$, these definitions are independent of this choice.
Indeed, they can be described in an analogous way to Proposition 
\ref{finmatbruhatprop}:
\bpr \label{matbruhatprop}
Let $\m,\m'\in M_{(b_i),n;(c_j),n'}$.
\ben
\item
$\ell(\m)=\sum_{i\in[1,n],j}m_{i,j}m_{\leq i,\geq j}
-\sum_{i\in[1,n],j} \binom{m_{i,j}+1}{2}$.
\item $\m\leq\m'$ if and only if, for all $i\in\Z$,
\[ m_{\leq i,\geq j}\leq m_{\leq i,\geq j}',\
\forall j,\text{ with equality for $j\ll 0$.}\]
\item $\m\leq\m'$ if and only if, for all $j\in\Z$,
\[ m_{\geq i,\leq j}\leq m_{\geq i,\leq j}',\
\forall i,\text{ with equality for $i\ll 0$.}\]
\een
\epr
\bpf
The proof is mostly identical to that of Proposition \ref{finmatbruhatprop},
using \eqref{lengtheqn} and Proposition \ref{bruhatprop} 
instead of \eqref{finlengtheqn} and Proposition \ref{finbruhatprop}.
In the proof of (3), the argument using \eqref{rewriteeqn} no longer 
makes sense, but the argument using transposes does.
\epf

As in \S2, we define $P_{\m,\m'}=P_{w_\m,w_{\m'}}$ for 
$\m,\m'\in M_{(b_i),n;(c_j),n'}$. From (2) of Proposition \ref{matbruhatprop}
it is clear that each interval $[\m,\m']$ in the poset
$M_{(b_i),n;(c_j),n'}$ is finite, so the inverse matrix
$(P_{\m,\m'}^{\langle -1\rangle})_{\m,\m'\in M_{(b_i),n;(c_j),n'}}$
of $(P_{\m,\m'})_{\m,\m'\in M_{(b_i),n;(c_j),n'}}$ is well defined.

The matrix definition of cancellability is identical to the finite case.
\bdf
If $\m\leq\m'$ in $M_{(b_i),n;(c_j),n'}$, we say that $(i,j)\in\Z\times\Z$
is \emph{cancellable} for $[\m,\m']$ if
\ben
\item $m_{i,j}\geq 1$.
\item $m_{\leq i-1,\geq j}=m_{\leq i-1,\geq j}'$, or equivalently
$m_{\geq i,\leq j+1}=m_{\geq i,\leq j+1}'$.
\item $m_{\leq i,\geq j+1}=m_{\leq i,\geq j+1}'$, or equivalently
$m_{\geq i+1,\leq j}=m_{\geq i+1,\leq j}'$.
\een
These equivalences follow from \eqref{affrewriteeqn}, bearing in mind
(2) of Proposition \ref{matbruhatprop}.
Clearly $(i,j)$ is cancellable iff $(i+n,j+n')$ is.
\edf
\bpr \label{readyprop}
Suppose that $\m\leq\m'$ in $M_{(b_i),n;(c_j),n'}$ and $(i,j)$ is 
cancellable for $[\m,\m']$. Let $\mathbf{e}$ be the matrix with
$e_{i+kn,j+kn'}=1$ for all $k$, all other entries zero.
\ben
\item For any $\m^1\in[\m,\m']$,
\ben
\item $m_{i,j}^1\geq m_{i,j}$,
\item $m_{\leq i-1,\geq j}^1=m_{\leq i-1,\geq j}$, and
\item $m_{\leq i,\geq j+1}^1=m_{\leq i,\geq j+1}$.
\een
Hence $(i,j)$ is cancellable for any sub-interval of $[\m,\m']$.
\item The map $\m^1\mapsto\m^1-\mathbf{e}$ is an isomorphism of posets
between $[\m,\m']$ and $[\m-\mathbf{e},\m'-\mathbf{e}]$, which reduces
all lengths by the same amount.
\item For any $\m^1,\m^2\in[\m,\m']$, 
$P_{\m^1,\m^2}=P_{\m^1-\mathbf{e},\m^2-\mathbf{e}}$.
\item For any $\m^1,\m^2\in[\m,\m']$, 
$P_{\m^1,\m^2}^{\langle -1\rangle}=
P_{\m^1-\mathbf{e},\m^2-\mathbf{e}}^{\langle -1\rangle}$.
\een
\epr
\bpf
Completely analogous to the proof of Proposition \ref{finreadyprop},
using the analogue of Proposition \ref{finmatcancprop}.
\epf
\section{Nilpotent orbits of the cyclic quiver}
We now return to the set-up of the latter part of \S1, so
$V$ is a $d$-dimensional 
$\Z/n\Z$-graded vector space, and $d_{i}=\dim V_{\bar{i}}$ for all 
$i\in\Z$.
We saw in \S1 that the $G_V$-orbits in $\CalN_V$ are in bijection with the
set $M_{(d_i),n}$ of multisegments (in the modulo $n$ sense) such that
each congruence class $\bar{i}$ occurs $d_i$ times among the elements of the
segments.
As in \S3, we will identify each $\m\in M_{(d_i),n}$ with a matrix
$(m_{i,j})$, this time in $M_{(d_i),n;(d_j),n}$; the definition
of $m_{i,j}$ is exactly the same as \eqref{conditioneqn}.
The resulting subset of $M_{(d_i),n;(d_j),n}$ is described as follows.
\bpr \label{idealprop}
Let $M_{(d_i),n}''=\{\m\in M_{(d_i),n;(d_j),n}\,|\,m_{i,j}=0,\,
\forall j<i-1\}$.
\ben
\item $M_{(d_i),n}''$ is a lower ideal of the poset $M_{(d_i),n;(d_j),n}$.
\item For all $\m\in M_{(d_i),n}''$, there is some $f(\m)\in\Z$ such that
\[ m_{\leq i,\geq j}+f(\m)=d_j+\cdots+d_i,\,\forall i\geq j. \]
\item For $\m,\m'$ in $M_{(d_i),n}''$, $\m\leq\m'$ if and only if 
$f(\m)=f(\m')$ and $m_{\leq i,\geq j}\leq m_{\leq i,\geq j}'$ for all
$i<j$.
\item $M_{(d_i),n}'=\{\m\in M_{(d_i),n}''\,|\,f(\m)=0\}$ is a lower ideal
of the poset $M_{(d_i),n;(d_j),n}$.
\item If $\m\in M_{(d_i),n}'$, then $m_{i,i-1}=m_{\leq i-1,\geq i}$ for all 
$i$.
\item $M_{(d_i),n}'=M_{(d_i),n}$.
\een
\epr
\bpf
As in the finite case, (1) is immediate from (3) of
Proposition \ref{matbruhatprop}. (2) comes from the fact that for 
$i\geq j$,
\[ m_{\leq i,\geq j}=d_j+\cdots+d_{i-1}+m_{\leq i,\geq i}
= d_{j+1}+\cdots+d_i+m_{\leq j,\geq j}. \]
Using this, (3) comes from
(2) of Proposition \ref{matbruhatprop}, and (4) is an immediate consequence
of (1) and (3). (5) is proved in the same way as (4) of Proposition
\ref{finidealprop}. From (5) and the $i=j$ case of (2) it follows that 
every matrix in $M_{(d_i),n}'$ arises from
a multisegment in $M_{(d_i),n}$, whence (6).
\epf

We can now state Lusztig's affine analogue of Theorem \ref{zelevinskythm}:
\bth \label{lusztigthm}
Let $\m,\m'\in M_{(d_i),n}$.
\ben
\item $\dim\CalO_\m=\ell(\m)-\sum_{i\in [1,n]}\binom{d_i}{2}$.
\item $\CalO_\m\subseteq
\overline{\CalO_{\m'}}\Leftrightarrow \m\leq \m'$.
\item $\CalH^i IC(\overline{\CalO_{\m'}})=0$ for $i$ odd. 
\item $IC_{\m,\m'}=P_{\m,\m'}$.
\item $IC_{\m,\m'}^{\langle -1\rangle}=P_{\m,\m'}^{\langle -1\rangle}$.
\een
\eth
\bpf
As with Theorem \ref{zelevinskythm}, (5) follows from (4) because
$M_{(d_i),n}$ is a lower ideal of $M_{(d_i),n;(d_j),n}$. Parts (1)--(4)
were proved by Lusztig in \cite[\S11]{quivers1}, but
since the conventions there are slightly different,
a sketch of a proof along the lines of the above proof of
Theorem \ref{zelevinskythm}
may be helpful.

Form $\CalV=\C(\!(t)\!)\otimes_\C V$, and consider lattices
(free $\C[\![t]\!]$-submodules of rank $d$) in $\CalV$.
Define $\Bh_{(d_i),n}$ to be the set of collections of 
lattices $(\CalM_i)_{i\in\Z}$
such that for all $i\in\Z$:
\ben
\item $\CalM_{i-1}\subset\CalM_{i}$,
\item $\dim_\C \CalM_{i}/\CalM_{i-1}=d_i$, and
\item $\CalM_{i-n}=t\CalM_i$.
\een 
It is well known that $\Bh_{(d_i),n}$ has the structure of an 
increasing union of projective varieties. We define a base-point
$(\CalL_i)$ in $\Bh_{(d_i),n}$ as follows.
For any $i\in\Z$, let $V_i$ denote $t^k V_{\bar{i}}$ where $k$
is defined by $i+kn\in\{1,\cdots,n\}$. Define
\[ \CalL_i=\widehat{\bigoplus_{j\leq i}}\, V_j,\ \forall i\in\Z, \]
where $\widehat{\bigoplus}$ denotes completed direct sum.
Relative to this base-point, $\Bh_{(d_i),n}$ decomposes into
affine Schubert cells $\Bh_\m$ for $\m\in M_{(d_i),n;(d_j),n}$.
Explicitly,
$\Bh_\m$ consists of those $(\CalM_i)$ such that for all $i,j\in\Z$,
\[ \dim_\C \frac{\CalM_i\cap\CalL_j}{\CalM_i\cap\CalL_{j-1}+
\CalM_{i-1}\cap\CalL_j} =m_{i,j}. \]
The analogues of (1)--(4) for these affine Schubert cells
(for all of $M_{(d_i),n;(d_j),n}$) are well known.
Let $\Bh_{(d_i),n}'$ be the closed subvariety of $\Bh_{(d_i),n}$
defined by requiring
$\CalM_i\supseteq\CalL_{i-1}$, $\dim_\C \CalM_i/\CalL_{i-1}=d_i$;
from the description of $M_{(d_i),n}$ as $M_{(d_i),n}'$,
it is easy to see that
$\Bh_{(d_i),n}'=\bigcup_{\m\in M_{(d_i),n}}\Bh_\m$.

Now we define a morphism $\CalN_V\to \Bh_{(d_i),n}':
\varphi\mapsto(\CalM_i(\varphi))$ by the rule
\[ \CalM_i(\varphi)=\CalL_{i-1}\oplus\{v+\varphi(v)+
\varphi^2(v)+\cdots\,|\,v\in V_i\}. \]
(Since $\varphi$ is nilpotent, this sum is actually finite.)
An easy check shows that this morphism maps
$\CalO_\m$ into $\Bh_\m$ for all $\m\in M_{(d_i),n}$.
All that remains is to verify that it gives an isomorphism between
$\CalN_V$ and the open subvariety of $\Bh_{(d_i),n}'$
defined by requiring
\[ \CalM_i\cap\bigoplus_{i'>i}V_{i'}=0,\ \forall i\in\Z. \]
The ``dual'' statement to this is what is proved in \cite[\S11]{quivers1}.
\epf
Note that in contrast to the situation in 
\S3, the poset $M_{(d_i),n}$ may have more than one maximal element.

We now come to the affine analogue of Theorem \ref{finmainthm},
a generalization of Theorem \ref{lusztigthm}.
\bth \label{mainthm}
Let $b_i,c_j\in\N$ be such that 
\[ b_{i+n}=b_i,\ c_{j+n}=c_j,\text{ and }d_i-b_i=d_{i-1}-c_{i-1},\,
\forall i,j\in\Z. \]
Define a $(\Z\times\Z)$-matrix $\mathbf{a}$ by
\[ a_{i,j}=\left\{\begin{array}{cl}
d_i-b_i,&\text{ if $j=i-1$,}\\
0,&\text{ otherwise.}
\end{array}\right. \]
Let $M_{(d_i),n}^{(b_i);(c_j)}=\{\m\in M_{(d_i),n}\,|\,m_{i,i-1}\geq d_i-b_i,\,
\forall i\in\Z\}$.
\ben
\item $M_{(d_i),n}^{(b_i);(c_j)}$ is an upper ideal of $M_{(d_i),n}$.
\item The map $\m\mapsto\m-\mathbf{a}$ is an isomorphism of posets between
$M_{(d_i),n}^{(b_i);(c_j)}$ and $\{\widetilde{\m}\in M_{(b_i),n;(c_j),n}\,|\,
\widetilde{m}_{i,j}=0,\,\forall j< i-1\}$.
\item For any $\m',\m''\in M_{(d_i),n}^{(b_i);(c_j)}$, 
$IC_{\m',\m''}=P_{\m'-\mathbf{a},\m''-\mathbf{a}}$.
\item For any $\m',\m''\in M_{(d_i),n}^{(b_i);(c_j)}$, 
$IC_{\m',\m''}^{\langle -1\rangle}=
P_{\m'-\mathbf{a},\m''-\mathbf{a}}^{\langle -1\rangle}$.
\een
\eth
\bpf
Completely analogous to the proof of Theorem \ref{finmainthm},
using Proposition \ref{matbruhatprop}, Proposition \ref{readyprop},
and Theorem \ref{lusztigthm} in place of Proposition \ref{finmatbruhatprop}, 
Proposition \ref{finreadyprop}, and Theorem \ref{zelevinskythm}.
\epf
\bcr \label{maincor}
For $\m\in M_{(d_i),n}$, define a $(\Z\times\Z)$-matrix $\m^-$ by
\[ m_{i,j}^-=\left\{\begin{array}{cl}
m_{i,j},&\text{ if $j=i-1$,}\\
0,&\text{ otherwise.}
\end{array}\right. \]
Let $\langle\m\rangle=\{\m'\in M_{(d_i),n}'\,|\,\m\leq\m'\}$, and
let $\m_1^{\mathrm{max}},\cdots,\m_t^{\mathrm{max}}$ be the 
maximal elements of $\langle\m\rangle$.
\ben
\item For all $\m'\in\langle\m\rangle$, 
$\m'-\m^-\in M_{(d_i-m_{i,i-1}),n;(d_j-m_{j+1,j}),n}$.
\item The map $\m'\mapsto\m'-\m^-$ is an isomorphism of posets between
$\langle\m\rangle$ and $\bigcup_{s=1}^t [\m-\m^-,\m_s^{\mathrm{max}}-\m^-]$.
\item For any $\m',\m''\in\langle\m\rangle$, 
$IC_{\m',\m''}=P_{\m'-\m^-,\m''-\m^-}$.
\item For any $\m',\m''\in\langle\m\rangle$, 
$IC_{\m',\m''}^{\langle -1\rangle}=P_{\m'-\m^-,\m''-\m^-}^{\langle -1\rangle}$.
\een
\ecr
\bpf
Apply Theorem \ref{mainthm} with
\[ b_i=d_i-m_{i,i-1},\ c_{i-1}=d_{i-1}-m_{i,i-1} \]
for all $i\in\Z$, and restrict to the upper ideal $\langle\m\rangle$
of $M_{(d_i),n}^{(b_i);(c_j)}$.
\epf
Note that the polynomials $P_{\m'-\m^-,\m''-\m^-}$ in (3) are
Kazhdan-Lusztig polynomials for $\widetilde{S_{k(\m)}}$, where
$k(\m)$ is the number of segments in $\m$, which is also the sum of the
entries in rows $1$ to $n$ of $\m-\m^-$.
As a corollary, we 
recover the main result of \cite{mytworow}:
\bcr \label{tworowcor}
If $\m\leq\m'$ in $M_{(d_i),n}$, and $k(\m)=2$,
then $IC_{\m,\m'}=1$.
\ecr
\bpf
In $\widetilde{S_2}$ all nonzero Kazhdan-Lusztig 
polynomials are $1$.
\epf

\bex
Let $d=6$, $n=3$, $d_1=d_2=d_3=2$. Let $\m\in M_{(d_i),3}$ be the
multisegment $[1,2]+[2,3]+[3,4]$. Then
$\langle\m\rangle$ has three maximal elements,
\[ \m_1^{\mathrm{max}}= [1,6],\
\m_2^{\mathrm{max}}= [2,7],\text{ and }
\m_3^{\mathrm{max}}= [3,8]. \]
Displaying only the
rows indexed by $1,2,3$, we have
\bes
\begin{split}
\m&=\begin{pmatrix}
\cdots&1&\mathbf{0}&1&0&0&0&0&\cdots\\
\cdots&0&1&0&1&0&0&0&\cdots\\
\cdots&0&0&1&0&1&0&0&\cdots
\end{pmatrix},\\
\m_1^{\mathrm{max}}&=\begin{pmatrix}
\cdots&1&\mathbf{0}&0&0&0&0&1&\cdots\\
\cdots&0&2&0&0&0&0&0&\cdots\\
\cdots&0&0&2&0&0&0&0&\cdots
\end{pmatrix},\\
\m-\m^-&=\begin{pmatrix}
\cdots&0&\mathbf{0}&1&0&0&0&0&\cdots\\
\cdots&0&0&0&1&0&0&0&\cdots\\
\cdots&0&0&0&0&1&0&0&\cdots
\end{pmatrix},\\
\m_1^{\mathrm{max}}-\m^-&=\begin{pmatrix}
\cdots&0&\mathbf{0}&0&0&0&0&1&\cdots\\
\cdots&0&1&0&0&0&0&0&\cdots\\
\cdots&0&0&1&0&0&0&0&\cdots
\end{pmatrix},
\end{split}
\ees
where the $\mathbf{0}$ is the $(1,1)$-entry.
Setting $B_i=\{i\}$ and $C_j=\{j-1\}$, so that $w_{\m-\m^-}$ is the
identity of $\widetilde{S_3}$, we find that $w_{\m_1^{\mathrm{max}}-\m^-}$
has window $(5,0,1)$, and is therefore $s_1 s_0 s_2 s_1$.
Similarly $w_{\m_2^{\mathrm{max}}-\m^-}=s_2 s_1 s_0 s_2$ and
$w_{\m_3^{\mathrm{max}}-\m^-}=s_0 s_2 s_1 s_0$.
So $\m'\mapsto w_{\m'-\m^-}$ is an isomorphism 
between $\langle\m\rangle$ and
$[1,s_1 s_0 s_2 s_1]\cup [1,s_2 s_1 s_0 s_2] \cup [1,s_0 s_2 s_1 s_0]$.
Moreover,
\[ IC_{\m,\m_1^{\mathrm{max}}}=P_{1,s_1 s_0 s_2 s_1}=q+1, \]
and similarly $IC_{\m,\m_2^{\mathrm{max}}}=IC_{\m,\m_3^{\mathrm{max}}}=q+1$,
while
\[ IC_{\m,\m'}^{\langle -1\rangle}=
P_{1,w_{\m'-\m^-}}^{\langle -1\rangle}=\varepsilon(w_{\m'-\m^-}) \]
for all $\m'\in\langle\m\rangle$ (an example of \eqref{affsigninveqn}).
\eex

Finally, we must prove Theorem \ref{introaffmainthm}.
We have elements
$\lambda,\mu\in \widetilde{D_k}$; define $\lambda_s$ and $\mu_s$ for all
$s\in\Z$ by the rule
\[ \lambda_{s+k}=\lambda_s+k-n,\ \mu_{s+k}=\mu_s+k-n. \]
Then $\lambda_s-s\geq\lambda_{s+1}-(s+1)$ for all $s\in\Z$,
and similarly for $\mu$; also,
\[ (w\cdot\mu)_s-s=\mu_{w^{-1}(s)}-w^{-1}(s),\ \forall w\in\widetilde{S_k},\,
s\in [1,k]. \]
Define
\[ B_i=\{-s\,|\,\mu_s-s+1=i\},\ C_j=\{-s\,|\,\lambda_s-s=j\}. \]
Then the subgroups $S_{(B_i)}$ and $S_{(C_j)}$ of $\widetilde{S_k}$
are exactly the images of $\widetilde{W_\mu}$ and $\widetilde{W_\lambda}$
under the automorphism $\tau:\widehat{S_k}\to\widehat{S_k}$ defined by
$\tau(w)(i)=-w(-i)$. 
The map $\psi:\widehat{S_k}\to M_{(b_i),n;(c_j),n}$
as defined in \S4 satisfies
\beq \label{psitaueqn}
\psi(\tau(w))_{i,j}=|\{s\in\Z\,|\,\mu_s-s+1=i,
\lambda_{w(s)}-w(s)=j\}|.
\eeq
So $\widetilde{S_k}[\lambda,\mu]=
\{w\in \widetilde{S_k}\,|\,\psi(\tau(w))_{i,j}=0,\,\forall
j<i-1\}$, which shows that it is indeed a lower ideal of $\widetilde{S_k}$.
Moreover, $w\mapsto \psi(\tau(w))$ gives an 
isomorphism of posets between
$\widetilde{S_k}[\lambda,\mu]^\circ$ and
\[ \{\widetilde{\m}\in M_{(b_i),n;(c_j),n}\,|\, 
\widetilde{\m}\geq \psi(1)\text{ and }\widetilde{m}_{i,j}=0,\,\forall
j<i-1\}, \]
and the polynomials attached to these posets correspond, since
\beq 
P_{w,w'}=P_{\tau(w),\tau(w')}
=P_{\psi(\tau(w)),\psi(\tau(w'))}
\eeq
for all $w,w'\in \widetilde{S_k}[\lambda,\mu]^\circ$. Now make the
assumption of Theorem 
\ref{introaffmainthm}, that $\lambda/\mu\in M_{(d_i),n}$; it follows
immediately that $(b_i)$ and $(c_j)$ satisfy the conditions of
Theorem \ref{mainthm}. By \eqref{psitaueqn},
for all $w\in \widetilde{S_k}[\lambda,\mu]^\circ$,
the multisegment $\lambda/(w\cdot\mu)$ when viewed as a matrix has
the same diagonal and above-diagonal entries as $\psi(\tau(w))$;
hence 
$\lambda/(w\cdot\mu)=\psi(\tau(w))+\mathbf{a}$
where $\mathbf{a}$ is as in Theorem \ref{mainthm}.
Thus Theorem \ref{introaffmainthm} follows from Theorem \ref{mainthm}.

\end{document}